\documentstyle[11pt]{article}

\input amssym 
\catcode`\@=11\relax
\newwrite\@unused
\def\typeout#1{{\let\protect\string\immediate\write\@unused{#1}}}

\def\psglobal#1{
\immediate\special{ps: plotfile #1 }}
\def\psfiginit{\typeout{psfiginit}
\immediate\psglobal{figtex.pro}%
\special{ps:: /TeXMagnification {\the\mag} def}
}

\def\@nnil{\@nil}
\def\@empty{}
\def\@psdonoop#1\@@#2#3{}
\def\@psdo#1:=#2\do#3{\edef\@psdotmp{#2}\ifx\@psdotmp\@empty \else
    \expandafter\@psdoloop#2,\@nil,\@nil\@@#1{#3}\fi}
\def\@psdoloop#1,#2,#3\@@#4#5{\def#4{#1}\ifx #4\@nnil \else
       #5\def#4{#2}\ifx #4\@nnil \else#5\@ipsdoloop #3\@@#4{#5}\fi\fi}
\def\@ipsdoloop#1,#2\@@#3#4{\def#3{#1}\ifx #3\@nnil
       \let\@nextwhile=\@psdonoop \else
      #4\relax\let\@nextwhile=\@ipsdoloop\fi\@nextwhile#2\@@#3{#4}}
\def\@tpsdo#1:=#2\do#3{\xdef\@psdotmp{#2}\ifx\@psdotmp\@empty \else
    \@tpsdoloop#2\@nil\@nil\@@#1{#3}\fi}
\def\@tpsdoloop#1#2\@@#3#4{\def#3{#1}\ifx #3\@nnil
       \let\@nextwhile=\@psdonoop \else
      #4\relax\let\@nextwhile=\@tpsdoloop\fi\@nextwhile#2\@@#3{#4}}
\def\psdraft{
	\def\@psdraft{0}
	\def\@psdraftspecial{100}
}
\def\psdraftspecial{
	\def\@psdraft{0}
	\def\@psdraftspecial{0}
}
\def\psfull{
	\def\@psdraft{100}
}
\psfull

\newif\if@prologfile
\newif\if@postlogfile
\newif\if@bbllx
\newif\if@bblly
\newif\if@bburx
\newif\if@bbury
\newif\if@height
\newif\if@width
\newif\if@rheight
\newif\if@rwidth
\newif\if@clip
\newif\if@right
\newif\if@left
\newif\if@toplines
\newif\if@box
\newif\if@caption
\newif\if@surround
\newif\if@captionwidth
\newif\if@captionwrite
\newif\if@captionopen
\def\@p@@sclip#1{\@cliptrue}
\def\@p@@sfile#1{
		\def\@p@sfile{#1}
}
\def\@p@@sfigure#1{
		\def\@p@sfile{#1}
}
\def\@p@sfake{\hbox to 0pt{\hss Whatever\hss}}
\def\@p@@sbbllx#1{
		\@bbllxtrue
		\@d@mscratch=#1
		\edef\@p@sbbllx{\number\@d@mscratch}
}
\def\@p@@sbblly#1{
		\@bbllytrue
		\@d@mscratch=#1
		\edef\@p@sbblly{\number\@d@mscratch}
}
\def\@p@@sbburx#1{
		\@bburxtrue
		\@d@mscratch=#1
		\edef\@p@sbburx{\number\@d@mscratch}
}
\def\@p@@sbbury#1{
		\@bburytrue
		\@d@mscratch=#1
		\edef\@p@sbbury{\number\@d@mscratch}
}
\def\@p@@sheight#1{
		\@heighttrue
		\@d@mscratch=#1
   		\edef\@p@sheight{\number\@d@mscratch}
}
\def\@p@@swidth#1{
		\@widthtrue
		\@d@mscratch=#1
		\edef\@p@swidth{\number\@d@mscratch}
}
\def\@p@@srheight#1{
		\@rheighttrue
		\@d@mscratch=#1
		\edef\@p@srheight{\number\@d@mscratch}
}
\def\@p@@srwidth#1{
		\@rwidthtrue
		\@d@mscratch=#1
		\edef\@p@srwidth{\number\@d@mscratch}
}
\def\@p@@sright#1{\@righttrue \@surroundtrue}
\def\@p@@sleft#1{\@lefttrue \@surroundtrue}
\def\@p@@sextraheight#1{\@d@mextraheight=#1}
\def\@p@@sbox#1{\@boxtrue}
\def\@p@@scaption#1{\@captiontrue}
\def\@p@@stoplines#1{
		\@toplinestrue
		\@c@ttoplines=#1
}
\def\@p@@scaptionwidth#1{
		\@captionwidthtrue
	  	\@d@mcaptionwidth=#1
}
\def\@p@@scaptionwrite#1{
		\global\@captionwritetrue
		\global\@w@rname=\expandafter{\jobname_captions.tex}
		\typeout{Captions are written to \the\@w@rname}
}

\def\@p@@sprolog#1{\@prologfiletrue\def\@prologfileval{#1}}
\def\@p@@spostlog#1{\@postlogfiletrue\def\@postlogfileval{#1}}
\def\@cs@name#1{\csname #1\endcsname}
\def\@setparms#1=#2,{\@cs@name{@p@@s#1}{#2}}
\def\ps@init@parms{
		\@bbllxfalse \@bbllyfalse
		\@bburxfalse \@bburyfalse
		\@heightfalse \@widthfalse
		\@rheightfalse \@rwidthfalse
		\def\@p@sbbllx{}\def\@p@sbblly{}
		\def\@p@sbburx{}\def\@p@sbbury{}
		\def\@p@sheight{}\def\@p@swidth{}
		\def\@p@srheight{}\def\@p@srwidth{}
		\def\@p@sfile{}
		\def\@p@scost{10}
		\def\@sc{}
		\@prologfilefalse
		\@postlogfilefalse
		\@clipfalse
		\@rightfalse \@leftfalse
		\@boxfalse \@captionfalse
		\@toplinesfalse \@surroundfalse
		\@d@mextraheight=0pt
 		\@c@ttoplines=0
		\@pshape={} \def\@p@srheight@total{}
		\@captionwidthfalse \@d@mcaptionwidth=0pt
}
\def\parse@ps@parms#1{
	 	\@psdo\@psfiga:=#1\do
		   {\expandafter\@setparms\@psfiga,}}
\newif\ifno@bb
\newif\ifnot@eof
\newread\ps@stream
\newtoks\@linetok
\def\bb@missing{
	\typeout{psfig: searching \@p@sfile \space  for bounding box}
	\openin\ps@stream=\@p@sfile
	\no@bbtrue
	\not@eoftrue
	\catcode`\%=12
	\loop
		\read\ps@stream to \line@in
		\global\@linetok=\expandafter{\line@in}
		\ifeof\ps@stream \not@eoffalse \fi
		\@bbtest{\@linetok}
		\if@bbmatch\not@eoffalse\expandafter\bb@cull\the\@linetok\fi
	\ifnot@eof \repeat
	\catcode`\%=14
}	
\catcode`\%=12
\newif\if@bbmatch
\def\@bbtest#1{\expandafter\@a@\the#1
\long\def\@a@#1
     \ifx\@bbtest#2\@bbmatchfalse\else\@bbmatchtrue\fi}
\long\def\bb@cull#1 #2 #3 #4 #5 {
	\@d@mscratch=#2 bp\edef\@p@sbbllx{\number\@d@mscratch}
	\@d@mscratch=#3 bp\edef\@p@sbblly{\number\@d@mscratch}
	\@d@mscratch=#4 bp\edef\@p@sbburx{\number\@d@mscratch}
	\@d@mscratch=#5 bp\edef\@p@sbbury{\number\@d@mscratch}
	\no@bbfalse
}
\catcode`\%=14
\def\compute@bb{
		\no@bbfalse
		\if@bbllx \else \no@bbtrue \fi
		\if@bblly \else \no@bbtrue \fi
		\if@bburx \else \no@bbtrue \fi
		\if@bbury \else \no@bbtrue \fi
		\ifno@bb \bb@missing \fi
		\ifno@bb \typeout{FATAL ERROR: no bb supplied or found}
			\no-bb-error
		\fi
		\count203=\@p@sbburx
		\count204=\@p@sbbury
		\advance\count203 by -\@p@sbbllx
		\advance\count204 by -\@p@sbblly
		\edef\@bbw{\number\count203}
		\edef\@bbh{\number\count204}
}
\def\in@hundreds#1#2#3{\count240=#2 \count241=#3
		     \count100=\count240	
		     \divide\count100 by \count241
		     \count101=\count100
		     \multiply\count101 by \count241
		     \advance\count240 by -\count101
		     \multiply\count240 by 10
		     \count101=\count240	
		     \divide\count101 by \count241
		     \count102=\count101
		     \multiply\count102 by \count241
		     \advance\count240 by -\count102
		     \multiply\count240 by 10
		     \count102=\count240	
		     \divide\count102 by \count241
		     \count200=#1\count205=0
		     \count201=\count200
			\multiply\count201 by \count100
		     	\advance\count205 by \count201
		     \count201=\count200
			\divide\count201 by 10
		     	\multiply\count201 by \count101
			\advance\count205 by \count201
		     \count201=\count200
			\divide\count201 by 100
			\multiply\count201 by \count102
			\advance\count205 by \count201
		     \edef\@result{\number\count205}
}
\def\compute@wfromh{
		\in@hundreds{\@p@sheight}{\@bbw}{\@bbh}
		\edef\@p@swidth{\@result}
}
\def\compute@hfromw{
		\in@hundreds{\@p@swidth}{\@bbh}{\@bbw}
		\edef\@p@sheight{\@result}
}
\def\compute@handw{
		\if@height
			\if@width
			\else
				\compute@wfromh
			\fi
		\else
			\if@width
				\compute@hfromw
			\else
				\edef\@p@sheight{\@bbh}
				\edef\@p@swidth{\@bbw}
			\fi
		\fi
}
\def\compute@resv{
		\if@rheight \else \edef\@p@srheight{\@p@sheight} \fi
		\if@rwidth \else \edef\@p@srwidth{\@p@swidth} \fi
		\edef\@p@srheight@total{\@p@srheight}
}
\newtoks\@pshape
\def\@c@ttoplines{\count120}
\def\@c@theightcount{\count121}
\def\@c@tshapecount{\count122}
\newdimen\@d@mwidthshape
\newdimen\@d@mextraheight
\newdimen\@d@mscratch
\def\compute@parshape{
	\if@right
		\if@left
	   		\typeout{error: Can't have both left and right set}
			\@leftfalse
		\fi
	\fi
	\@d@mscratch=\@p@swidth truesp
	\divide \@d@mscratch by 19
	\multiply \@d@mscratch by 20
	\edef\@p@swidthdimen{\the\@d@mscratch}
	\@c@tshapecount=\@c@ttoplines
 	\@d@mscratch=\baselineskip
	\multiply \@d@mscratch by \@c@ttoplines
	\advance \@d@mscratch by .4\baselineskip
    	\edef\@p@stopdistance{\the\@d@mscratch }
	\@d@mscratch=\@p@sheight truesp
	\divide \@d@mscratch by 2
	\edef\@p@shalfboxheight{\the\@d@mscratch}
	\if@toplines
		\loop \@pshape=\expandafter{\the\@pshape 0pt \hsize}
		\advance\@c@ttoplines by -1
		\ifnum\@c@ttoplines>0 \repeat
	\fi
   	\@c@theightcount=\@p@srheight@total
	\advance \@c@theightcount by \@d@mextraheight
	\divide  \@c@theightcount by \baselineskip
	\advance \@c@theightcount by 1
    	\advance \@c@tshapecount by \@c@theightcount
	\advance \@c@theightcount by -1
	\@d@mwidthshape=\hsize
     	\advance \@d@mwidthshape by -\@p@swidthdimen
	\if@left
		\def\@moveshape{0pt}
		\ifnum\@c@theightcount>0
		      	\loop
			\@pshape=%
			\expandafter{\the\@pshape %
					\@p@swidthdimen \@d@mwidthshape}
			\advance \@c@theightcount by -1
			\ifnum\@c@theightcount>0 \repeat
		\else
			\advance \@c@tshapecount by 1
		\fi
	\fi
	\if@right
		\@d@mscratch=\hsize
		\advance \@d@mscratch by -\@p@swidth truesp
		\edef\@moveshape{\@d@mscratch}
		\ifnum\@c@theightcount>0
			\loop
			\@pshape=\expandafter{\the\@pshape 0pt \@d@mwidthshape}
			\advance \@c@theightcount by -1
			\ifnum\@c@theightcount>0 \repeat
		\else
			\advance \@c@tshapecount by 1
		\fi
	\fi
	\ifnum \@p@srheight=0
		\@pshape={}
		\@c@tshapecount = 0
	\else
	 	\@pshape=\expandafter{\the\@pshape 0pt \hsize}
	\fi
}
\def\@p@ssetsurroundboxes{
	\global\parshape=\count122 \the\@pshape		
 	\moveright\@moveshape
	\vbox to 0pt\bgroup\hskip0pt\vskip\@p@stopdistance
}
\newtoks\@captiontok
\newbox\@b@xcaption
\newdimen\@d@mcaptionwidth
\newdimen\@d@mcaptionheight
\newwrite\@w@rcaption
\newtoks\@w@rname
\def\setcaption#1{\@captiontok={#1}}
\def\@set@caption{
	\setbox\@b@xcaption\vbox{\hsize\@d@mcaptionwidth
	\tolerance=9000 \baselineskip=11.4pt
	\noindent\relax\the\@captiontok}
	\if@captionwrite
		\if@captionopen
		\else
			\global\@captionopentrue
			\immediate\openout\@w@rcaption=\the\@w@rname
		\fi
		\immediate\write\@w@rcaption{\the\@captiontok}
		\immediate\write\@w@rcaption{}
	\fi
}
\def\compute@caption{
	\if@captionwidth
	\else
		\@d@mcaptionwidth = \@p@swidth truesp
		\divide \@d@mcaptionwidth by 20
		\multiply \@d@mcaptionwidth by 17
	\fi
	\@set@caption
	\@d@mcaptionheight=\ht\@b@xcaption
	\if@rheight
	\else
		\count100=\@p@srheight
	   	\advance \count100 by \@d@mcaptionheight
	   	\advance \count100 by \bigskipamount
	   	\advance \count100 by \medskipamount
	   	\edef\@p@srheight@total{\number\count100}
	\fi
}
\newif\if@alreadyjtem \@alreadyjtemfalse
\def\newpar{\hfil\vadjust{\vskip\parskip}%
	{\count100=\parskip
	\count101=\baselineskip
	\divide\count101 by 10  \multiply\count101 by 3
	\advance \count100 by \count101
	\divide\count100 by \baselineskip
	\advance\count100 by \prevgraf
	\global\prevgraf=\count100}%
	\break\if@alreadyjtem\else\indent\fi%
}
\let\sav@par=\par
\def\jtem#1{%
    	\if@alreadyjtem\else\bgroup\fi
	\def\par{\sav@par\egroup\sav@par}
	\if@alreadyjtem\else\leftskip\parindent\fi
	\@alreadyjtemtrue
	\noindent\hskip0pt
	\llap{#1\ }\ignorespaces
}
\def\compute@sizes{%
	\compute@bb
	\compute@handw
  	\compute@resv
	\if@caption
		\compute@caption
	\fi
	\if@surround
		\compute@parshape
	\fi
}
\def\@p@sdospecials{
	\ifnum\@p@scost<\@psdraft
	       	\typeout{psfig: including \@p@sfile \space }
	\fi
	\special{ps::[begin] 	\@p@swidth \space \@p@sheight \space
			\@p@sbbllx \space \@p@sbblly \space
			\@p@sbburx \space \@p@sbbury \space
			startTexFig \space }
	\ifnum\@p@scost<\@psdraft
		\if@clip
			\typeout{(clip)}
			\special{ps:: \@p@sbbllx \space \@p@sbblly \space
				\@p@sbburx \space \@p@sbbury \space
			    	doclip \space }
		\fi
	\fi
	\if@box
		\typeout{(box)}
  		\special{ps:: \@p@sbbllx \space \@p@sbblly \space
			\@p@sbburx \space \@p@sbbury \space
		    	dobox \space }
	\fi
	\ifnum\@p@scost<\@psdraft
		\if@prologfile
	    		\special{ps: plotfile \@prologfileval \space }
		\fi
		\special{ps: plotfile \@p@sfile \space }
    		\if@postlogfile
			\special{ps: plotfile \@postlogfileval \space }
		\fi
	\fi
	\special{ps::[end] endTexFig \space }
}
\newif\if@putinvbox

\def\psfig#1{{%
	\ifhmode%
		\vbox\bgroup
		\@putinvboxtrue
	\else
		\@putinvboxfalse
	\fi
       	\ps@init@parms
	\parse@ps@parms{#1}
       	\compute@sizes
	\if@surround
		\psfig@for@surround
	\else
		\psfig@for@regular
	\fi
	\if@putinvbox
       		\egroup
	\fi
}}
\def\psfig@for@surround{%
	\@p@ssetsurroundboxes
	\ifnum\@p@scost<\@psdraft
		\@p@sdospecials
		\vbox to \@p@srheight true sp{\vss}
       	\else
		\if@box
			\@p@sdospecials
		\fi
		\vbox to \@p@srheight true sp{
			\vskip\@p@shalfboxheight
			\hbox to \@p@srwidth true sp{
				\hss
				\ifnum\@p@scost<\@psdraftspecial
					\@p@sfile
				\else
					\@p@sfake
				\fi
      				\hss
			}
		\vss
		}
	\fi
	\if@caption
		\medskip
		\hbox to \@p@srwidth true sp{
			\hss
			\box\@b@xcaption
			\hss
		}
 		\medskip
	\fi
	\vss\egroup
	\vskip-\parskip
}

\def\psfig@for@regular{%
	\if@putinvbox
	\else
		\vskip\parskip
	\fi
	\ifnum\@p@scost<\@psdraft
		\@p@sdospecials
		\vbox to \@p@srheight true sp{%
			\hbox to \@p@srwidth true sp{
			\hfil
			}
		\vfil
		}
       	\else
		\if@box
			\@p@sdospecials
		\fi
	    	\vbox to \@p@srheight true sp{
			\vss
			\hbox to \@p@srwidth true sp{
				\hss
				\ifnum\@p@scost<\@psdraftspecial
					\@p@sfile
				\else
					\@p@sfake
				\fi
				\hss
			}
		    	\vss
		}
	\fi
	\if@caption
		\medskip
		\hbox to \@p@srwidth true sp{
			\hss
			\box\@b@xcaption
			\hss
		}
		\bigskip
	\fi
	\if@putinvbox
	\else
		\vskip-\parskip
	\fi
}
\catcode`\@=12\relax
  
\psfiginit

\font\tinybbfont=msbm6
\font\scriptsizebbfont=msbm7 scaled \magstep 1
 1
\font\footnotesizebbfont=msbm9 scaled \magstep 0
\font\smallbbfont=msbm7 scaled \magstep 2
\font\bbfont=msbm9 scaled \magstep1  
\font\largebbfont=msbm10 scaled \magstep 1
 2
 3
 4

\def\tinyBbb#1{\hbox{\tinybbfont #1}}
\def\scriptsizeBbb#1{\hbox{\scriptsizebbfont #1}}

\def\footnotesizeBbb#1{\hbox{\footnotesizebbfont #1}}
\def\smallBbb#1{\hbox{\smallbbfont #1}}
\def\Bbb#1{\hbox{\bbfont #1}}
\def\largeBbb#1{\hbox{\largebbfont #1}}

\newcommand{\Ann}{\mbox{\it Ann}}

\newcommand{\CP}{\mbox{{\Bbb C}{\rm P}}}

\newcommand{\CPscriptsize}{\mbox{\scriptsize {\scriptsizeBbb C}{\rm P}}}
\newcommand{\CPsmall}{\mbox{\small {\smallBbb C}{\rm P}}}

\newcommand{\Diag}{\mbox{\it Diag}\,}

\newcommand{\Fl}{\mbox{\it Fl}\,}
\newcommand{\GL}{\mbox{\it GL}\,}
\newcommand{\Gr}{\mbox{\it Gr}}
\newcommand{\Hom}{\mbox{\it Hom}\,}

\newcommand{\PGL}{\mbox{\it PGL}\,}
\newcommand{\Pic}{\mbox{\rm Pic}\,}
 
\newcommand{\Proj}{\mbox{\rm Proj}\,}

\newcommand{\Quot}{\mbox{\it Quot}\,}

\newcommand{\SL}{\mbox{\it SL}}

\newcommand{\Span}{\mbox{\it Span}\,}
\newcommand{\Spec}{\mbox{\it Spec}\,}

\newcommand{\Sym}{\mbox{\it Sym}}

\newcommand{\Wt}{\mbox{\it Wt}}

\newcommand{\degree}{\mbox{\it deg}\;}
\newcommand{\dimm}{\mbox{\it dim}\,}
\newcommand{\ev}{\mbox{\it ev}\,}

\newcommand{\pr}{\mbox{\rm pr}}
\newcommand{\pt}{\mbox{\it pt}}
\newcommand{\rank}{\mbox{\it rank}\,}

\begin{document}

\enlargethispage{23cm}

\begin{titlepage}

$ $

\vspace{-2cm} 

\noindent\hspace{-1cm}
\parbox{6cm}{\small October 2001}\
   \hspace{6.5cm}\
   \parbox{5cm}{math.AG/0111256}

\vspace{1.5cm}

\centerline{\large\bf
The $S^1$ fixed points in Quot-schemes}
\vspace{1ex}
\centerline{\large\bf and mirror principle computations}

\vspace{1.5cm}

\centerline{\large Bong H.\ Lian$^{1,\,a}$,\hspace{2ex}
            Chien-Hao Liu$^{2,\,b}$,\hspace{2ex}
            Kefeng Liu$^{3,\,c}$,\hspace{2ex} and\hspace{2ex}
            Shing-Tung Yau$^{4,\,b}$}
\vspace{1.1em}
\hspace{-7em}
\begin{minipage}{18cm}
 \begin{center}
  \parbox[t]{4cm}{
   \centerline{\it $^a$Department of Mathematics}
   \centerline{\it Brandeis University}
   \centerline{\it Waltham, MA 02154} } \
  \hspace{1em} \ 
  \parbox[t]{6cm}{
   \centerline{\it $^b$Department of Mathematics}
   \centerline{\it Harvard University}
   \centerline{\it Cambridge, MA 02138}  } \
  \hspace{1em} \
  \parbox[t]{5cm}{
   \centerline{\it $^c$Department of Mathematics}
   \centerline{\it University of California at Los Angelas}
   \centerline{\it Los Angelas, CA 90095}  } \
 \end{center}
\end{minipage}

\vspace{3em}

\begin{quotation}
\centerline{\bf Abstract}
\vspace{0.3cm}
\baselineskip 12pt  
{\small
We describe the $S^1$-action on the Quot-scheme $\Quot({\cal E}^n)$
 associated to the trivial bundle
 ${\cal E}^n=\CPsmall^1\times{\smallBbb C}^n$.
 In particlular, the topology of the $S^1$-fixed-point components
 in $\Quot({\cal E}^n)$ and the $S^1$-weights of the normal bundle
 of these components are worked out.
Mirror Principle, as developed by three of the current authors in
 the series of work [L-L-Y1, I, II, III, IV], is a method for studying
 certain intersection numbers on a stable map moduli space.
As an application, in Mirror Principle III, Sec 5.4, an outline was
 given in the case of genus zero with target a flag manifold.
The results on $S^1$ fixed points in this paper are used here to do
 explicit Mirror Principle computations in the case of Grassmannian
 manifolds. In fact,  Mirror Principle computations involve only
 a certain distinguished subcollection of the $S^1$-fixed-point
 components. These components are identified and are labelled
 by Young tableaus. The $S^1$-equivariant Euler class $e_{S^1}$ of
 the normal bundle to these components is computed. A diagrammatic
 rule that allows one to write down $e_{S^1}$ directly
 from the Young tableau is given. From this, the aforementioned
 intersection numbers on the moduli space of stable maps can be worked
 out. Two examples are given to illustrate our method.
Using our method, the A-model for Calabi-Yau complete intersections 
 in a Grassmannian manifold can now also be computed explicitly.
This work is motivated by the intention to provide further details
 of mirror principle and to understand the relation of mirror principle
 to physical theory. Some related questions are listed for further study.  
} 
\end{quotation}

\bigskip

\baselineskip 12pt
{\footnotesize
 } 

\medskip

\noindent {\small
MSC number 2000$\,$: 14N35, 14D20, 55N91, 55R91, 81T30.
} 

\bigskip

\baselineskip 11pt
{\footnotesize
\noindent{\bf Acknowledgements.}
 We would like to thank
  Raoul Bott, Jiun-Cheng Chen, Joe Harris, Shinobo Hosono,
   Yi Hu, Deepee Khosla, Albrecht Klemm,
   Jun Li, Hui-Wen Lin, Chiu-Chu Liu, Alina Marian,
   Keiji Oguiso, Chin-Lung Wang, and Xiaowei Wang
 for valuable conversations, related group meeting talks,
  interesting questions, and guidance of literatures
  that strongly affect the work.
 C.-H.\ L would like to thank in addition 
  Orlando Alvarez and William Thurston for educations;
  J.H., Sean Keel, Mihnea Popa, and Ravi Vakil
 for the various courses on algebraic geometry;
  Dima Arinkin, Hung-Wen Chang, J.-C.C., Izzet Coskun, M.P.,
  and Jason Starr
  for the valuable discussions on related objects
  in algebraic geometry at various stages of the work;
  and Ling-Miao Chou for the tremendous moral support.
 The work is supported by DOE grant DE-FG02-88ER25065,
  NSF grants DMS-9803347, and DMS-0072158.
} 

\noindent
\underline{\hspace{20em}}

\noindent
{\footnotesize E-mail:
$^1$lian@brandeis.edu$\,$,\hspace{1ex}
$^2$chienliu@math.harvard.edu$\,$,\hspace{1ex}
$^3$liu@math.ucla.edu$\,$,\hspace{1ex}
$^4$yau@math.harvard.edu}

\end{titlepage}

\newpage
$ $

\vspace{-4em}  

\centerline{\sc  Quot-Schemes and Mirror Principle}

\vspace{2em}

\baselineskip 14pt  

\begin{flushleft}
{\Large\bf 0. Introduction and outline.}
\end{flushleft}

\begin{flushleft}
{\bf Introduction.}
\end{flushleft}
In Mirror Principle III, Sec.\ 5.4, in the series of work
 [L-L-Y1, I, II, III, IV] developed by three of the current authors,
 they outlined how Mirror Principle can be used to study
 certain intersection numbers on a stable map moduli space for flag manifolds.
 In this article, we carry out this computation explicitly in the case of
 Grassmannian manifolds.
This is our main motivation for studying 
the $S^1$-action on Quot-schemes. The latter is, of course,
of independent interests from the viewpoint of group actions on manifolds,
regardless of Mirror Principle. Two of our main results are
the topology of the $S^1$-fixed-point components
 in $\Quot({\cal E}^n)$ (Theorem 2.1.9), 
and the $S^1$-weights of the normal bundle
 to these components (Theorem 2.2.1).
Mirror Principle computations involve only a certain distinguished
 subcollection of the $S^1$-fixed-point components.
 These components are identified and are labelled by Young tableaus.
The $S^1$-equivariant Euler class $e_{S^1}$ of the normal bundle
 of each of these components is computed (Theorem 3.3.3).
 A diagrammatic rule that allows one to write down $e_{S^1}$ directly
 from the Young tableau is given.
From this, the intersection numbers
  of the moduli space of stable maps can be easily worked out.
Two sample calculations are given to illustrate the method (Sec.\ 4).
 The answers are self-consistent and is the same as the result
 computed via the method of Mirror Principle I in a special case.
Using our method, the A-model for Calabi-Yau complete intersections 
 in a Grassmannian manifold can now also be computed explicitly.
This work is motivated by the intention to provide further details
 of mirror principle and to understand the relation of mirror principle
 to physical theory.
Some related questions are listed in the end for further pursuit.

\bigskip

\bigskip

\begin{flushleft}
{\bf Outline.}
\end{flushleft}
{\small
\baselineskip 11pt  
\begin{quote}
 1. Essential backgrounds and notations for physicists.

 2. The $S^1$-action on Quot-schemes.
    \vspace{-1ex}
    \begin{quote}
     \hspace{-1.3em}
     2.1 The $S^1$-fixed-point components.

     \hspace{-1.3em}
     2.2 \ \parbox[t]{12cm}{
         The $S^1$-weight system of the tangent space
         of Quot-scheme at an $S^1$-fixed-point component. }

     \hspace{-1.3em}
     2.3 Combinatorics of the $S^1$-weight system
         and the multiplicity of $0$.
    \end{quote}

 \vspace{-.8ex}
 3. Mirror principle computation for Grassmannian manifolds.
    \vspace{-1ex}
    \begin{quote}
     \hspace{-1.3em}
     3.1 \ \parbox[t]{12cm}{
         The distinguished $S^1$-fixed-point components
         and the hyperplane-induced class.}

     \medskip
     \hspace{-1.3em}
     3.2 \ \parbox[t]{12cm}{
          The weight subspace decomposition of the normal bundle
          to the distinguished components.}

     \hspace{-1.3em}
     3.3 Structure of the induced bundle and
         the ${\smallBbb C}^{\times}$-equivariant Euler class.
    \end{quote}

 \vspace{-.8ex}
 4. Illustrations by two examples.
\end{quote}
} 

\newpage
\baselineskip 14pt  

\section{Essential backgrounds and notations for physicists.}

Essential backgrounds or their main references used in this article
 and notations for objects involved are collected in this section
 for the convenience of readers.

\bigskip

\noindent $\bullet$
{\bf Schemes, coherent sheaves, and Hilbert polynomials.}
(See
  Eisenbud-Harris [E-H], Hartshorne [Ha$\,$: Chapter II], and
  Friedman [Fri$\,$: Chapter 2];
  also Kempf [Ke] and Mumford [Mu2].)
Let $X$ be a projective variety with a fixed very ample line bundle
 ${\cal O}(1)$, then the Hilbert polynomial of coherent sheaves on
 $X$ are additive with respect to short exact sequences.
In other words, if
 $0\rightarrow {\cal F}^{\prime}\rightarrow {\cal F}
             \rightarrow {\cal F}^{\prime\prime}\rightarrow 0$
 is an exact sequence of coherent sheaves on $X$,
 then
 $P_{{\cal O}(1)}({\cal F})
       =P_{{\cal O}(1)}({\cal F}^{\prime})
         + P_{{\cal O}(1)}({\cal F}^{\prime\prime})\,$,
 where $P_{{\cal O}(1)}(\,\cdot\,)$ is the Hilbert polynomial of
 $\,\cdot\,$. (Cf.\ [Ha$\,$: Ex.\ III.5.1]).

\bigskip

\noindent $\bullet$
{\bf Coherent sheaves on a curve and their Hilbert polynomial.}
\begin{itemize}
 \item [-]
  A coherent sheaf ${\cal F}$ on a smooth curve $C$ fits into
   a {\it split exact sequence} of ${\cal O}_C$-modules:
   $0\rightarrow {\cal F}_{\rm torsion} \rightarrow {\cal F}
           \rightarrow  {\cal F}^{\vee\vee}\rightarrow 0$,
   where ${\cal F}_{\rm torsion}$ is the torsion subsheaf of
   ${\cal F}$ and the double dual ${\cal F}^{\vee\vee}$ of
   ${\cal F}$ is locally-free.
  In case $C$ is nodal, then an exact sequence from
   the normalization $\widetilde{C}$ of $C$ can be used to
   understand coherent sheaves on $C$ as well.

 \item [-]
  Let ${\cal F}$ be a coherent sheaf on a smooth curve $C$, then
  $$
   \degree{\cal F}\; =\; c_1({\cal F})\;
    =\; c_1({\cal F}_{\rm torion})\,+\,c_1({\cal F}^{\vee\vee})\;
    =\; \dimm_{\scriptsizeBbb C}\,\Gamma(C,{\cal F}_{\rm torsion})\,
         +\, c_1({\cal F}^{\vee\vee})\,.
  $$

 \item [-]
  Fix a very ample line bundle ${\cal O}(1)$ on $C$, let $k$ be
  the rank of ${\cal F}$ and $g$ be the arithmetic genus of $C$.
  Then the Hilbert polynomial for ${\cal F}$ is given by
  $$
   P({\cal F},t)\;
    =\; (k\,\degree C)\, t\,+\, \deg{\cal F} + k(1-g)\,.
  $$
  For $C=\CP^1$ with ${\cal O}_{\CPscriptsize^1}(1)\,$,
   this is $P({\cal F},t)=\,k\,t + (c_1({\cal F})+k)\,$.
  For ${\cal F}$ a torsion sheaf, $r=0$ and the polynomial becomes
   $P({\cal F},t)\,=\,c_1({\cal F})\,
    =\,\dimm_{\scriptsizeBbb C}\,\Gamma(C,{\cal F})\,$.
\end{itemize}
Cf.\  [H-L], [H-M], [Ke], and [LP].

\bigskip

\noindent $\bullet$ {\bf Quot-scheme.}
(See Huybrechts-Lehn [H-L$\,$: Chapter 2];
 also Grothendieck [Gr3], Koll\'{a}r [Kol$\,$: Sec.\ I.1], and
 Mumford [Mu1].)
Let
 $S$ be a projective variety $S$ with a fixed ample line bundle,
 and ${\cal F}$ be a coherent sheaf on $S$.
Then the Quot-scheme $\Quot_P({\cal E}^n)$ of Grothendieck is
 the fine moduli space that parameterizes the set of quotients
 ${\cal F}
   \rightarrow {\cal F}\!/\mbox{\raisebox{-.4ex}{${\cal V}$}}$
 with $P({\cal F}\!/\mbox{\raisebox{-.4ex}{${\cal V}$}},t)$
 a given polynomial $P=P(t)$.
 It is the scheme that represents the $\Quot$-functor of
 Grothendieck, cf [Gr3].

\bigskip

\noindent $\bullet$
{\bf Quot-scheme compactification of $\Hom(\CP^1, \Gr_r({\Bbb C}^n))$.}
(Cf.\ [Str].)
 Let
  $C=\CP^1$ with the very ample line bundle
    ${\cal O}_{\CPscriptsize^1}(1)$,
  ${\cal E}^n$ be a trivialized trivial bundle of rank $n$ over $C$,
  $\Gr_r({\Bbb C}^n)$ be the Grassmannian manifold that parameterizes
    $r$-planes in ${\Bbb C}^n$,
    and
  $\Hom(\CP^1, \Gr_r({\Bbb C}))$ be the space of morphisms from
    $\CP^1$ to $\Gr_r({\Bbb C}^n)$.
 Then an element $(f:\CP^1\rightarrow \Gr_r({\Bbb C}^n))$ in
  $\Hom(\CP^1,\Gr_r({\Bbb C}))$ determines a unique rank-$r$ subbundle
  ${\cal V}$ in ${\cal E}^n$, which corrsponds in turn to the element
  ${\cal E}^n
     \rightarrow {\cal E}^n\!/\mbox{\raisebox{-.4ex}{${\cal V}$}}$
  in $\Quot({\cal E}^n)$.
 This gives a natural embedding of $\Hom(\CP^1,\Gr_r({\Bbb C}^n))$
  in $\Quot({\cal E}^n)$.
 The component of $\Hom(\CP^1,\Gr_r({\Bbb C}^n))$ that contains
  degree-$d$ image curves in $\Gr_r({\Bbb C}^n)$ is embedded in
  $\Quot_P({\cal E}^n)$ with the Hilbert polynomial
  $P=P(t)=(n-r)t+d+(n-r)$.
 This gives a compactification of $\Hom(\CP^1, \Gr_r({\Bbb C}^n))$
  via Quot-schemes, other than the moduli space of stable maps.
 Recall also that $\Quot_P({\cal E}^n)$ is a smooth, irreducible,
  rational projective variety of dimension $dn+(n-r)r$,
 cf.\ [Ch] and [Kim].
 The $S^1$-action on $\CP^1$ induces $S^1$-actions on
  $\Hom(\CP^1,\Gr_r({\Bbb C}^n))$ and $\Quot({\cal E}^n)$
  respectively.
 The two actions coincide under the natural embedding of
  $\Hom(\CP^1, \Gr_r({\Bbb C}^n))$ in $\Quot({\cal E}^n)$.

\bigskip

\noindent $\bullet$
{\bf Mirror principle for Grassmannian manifolds.}
For the details of Mirror Principle, readers are referred to
 [L-L-Y1$\,$: I, II, III, IV]. Some survey is given in [L-L-Y2].
To avoid digressing too far away, here we shall take
 [L-L-Y1, III$\,$: Sec.\ 5.4] as our starting point and
 restrict to the case that the target manifold of stable maps is
 $X=\Gr_r({\Bbb C}^n)$.
Recall first the Pl\"{u}cker embedding 
$\tau:X=\Gr_r({\Bbb C}^n)\,\rightarrow\,Y=\CP^{\,n\choose r}$,
which induces an isomorphism between the divisor class groups
$\tau^{\ast}:A^1(Y)\stackrel{\sim}{\rightarrow} A^1(X)$.

Recall next the setup of Mirror Principle for $X=\Gr_r({\Bbb C}^n)$.
The geometric objects involved are contained in the following
 diagram$\,$:

{\small
$$
 \begin{array}{cccccccclcl}
  V  &  & U_d & & V_d & & {\cal U}_d        \\
  \downarrow   & & \downarrow   & & \downarrow   & & \downarrow \\
  X  & \stackrel{ev}{\longleftarrow}
     & M_{0,1}(d,X)    & \stackrel{\rho}{\longrightarrow}
     & M_{0,0}(d, X)   & \stackrel{\pi}{\longleftarrow}
     & M_d             & \stackrel{\varphi}{\longrightarrow}
     & W_d             & \stackrel{\psi}{\longleftarrow}
     & \Quot_{(d)}                                       \\
   & & & & & & \cup    & & \hspace{.6ex}\cup  & & \hspace{.6ex}\cup \\
   & & & & & & F_0     & \stackrel{ev^Y}{\longrightarrow}
     & Y_0\,(\supset X_0=X)   & \stackrel{g}{\longleftarrow}
     & E_0\;=\;\cup_s\,E_{0s}\,,\\
   & & & & & & \mbox{\scriptsize $|$}\wr  & & \mbox{\scriptsize $|$}\wr \\
   & & & & & & X  & \stackrel{\tau}{\longrightarrow} & Y
 \end{array}
$$
{\normalsize where}} 
\begin{itemize}
 \item [(1)] {\it Moduli spaces}$\,$:
  $M_{0,0}(d,X)$ is the moduli space of genus-$0$ stable maps of
   degree $d$ into $X$,
  $M_{0,1}(d,X)$ is the moduli space of genus-$0$, $1$-pointed
   stable maps of degree $d$ into $X$,
  $M_d=M_{0,0}(\CP^1\times X, (1,d))$,
  $W_d$ is the linear sigma model at degree $d$,
   which can be chosen to be the projective space 
   ${\Bbb P}(H^0(\CP^1,\,{\cal O}_{{\scriptsizeBbb C}{\rm P}^1}(d))
                                   \otimes \Lambda^r{\Bbb C}^n)$
   for $X=\Gr_r({\Bbb C}^n)$, and
  $\Quot_{(d)}=\Quot_P({\cal E}^n)$ with $P=P(t)=(n-r)t+d+(n-r)$;

 \item [(2)] {\it Group actions}$\,:$
  there are ${\Bbb C}^{\times}$-actions on
   $M_d$, $W_d$, and $\Quot_{(d)}$ respectively that are compactible
   with the morphisms among these moduli spaces;
  these ${\Bbb C}^{\times}$-actions induce $S^1$-actions on these
   moduli spaces by taking the subgroup
   $U(1)\subset{\Bbb C}^{\times}$;

 \item [(3)] {\it Morphisms}$\,$:
   $\ev$ is the evaluation map,
   $\rho$ is the forgetful map,
   $\pi$ is the contracting morphism,
   $\varphi$ is the collapsing morphism, and
   $\psi$ is an $S^1$-equivariant resolution of singularities
    of $\varphi(M_d)$, which will be discuessed in detail
    in Sec.\ 3.1;

 \item [(4)] {\it Bundles}$\,$:
   $V$ is a vector bundle over $X$,
   $V_d=\rho_!\ev^{\ast}V$, $U_d=\rho^{\ast}V_d$,
   and ${\cal U}_d=\pi^{\ast}V_d$;

 \item [(5)] {\it Special $S^1$-fixed-point locus}$\,$:
   $F_0\simeq M_{0,1}(d,X)$ is the special $S^1$-fixed-point
    component in $M_d$ that corresponds to gluing stable maps
    $(C^{\prime}, f^{\prime}, x^{\prime})$ to $\CP^1$ at
    $x^{\prime}\in C^{\prime}$ and $\infty\in \CP^1$,
   $Y_0$ is the special $S^1$-fixed-point component in $W_d$
    such that $\varphi^{-1}(Y_0)=F_0$, and 
   $E_0$ is the $S^1$-fixed-point locus in $\psi^{-1}(Y_0)$
   and is called the distinguished $S^1$-fixed-point locus
   or components in $\Quot_P({\cal E}^n)$.
\end{itemize}

Associated to each $(V,b)$, where $b$ is a multiplicative
characteristic class, is the Euler series
$A(t)\in A^{\ast}(X)(\alpha)[t]\,$:
$$
 \begin{array}{lllll}
  A(t)  & =  & A^{V,b}  & =
    & e^{-H\cdot t/\alpha}\,\sum_d\,A_d\,e^{d\cdot t}\,, \\[1.2ex]
  A_d   & =  & i_0^{\ast}\,b({\cal U}_d) & :=
    & \ev^X_{\ast}\,\left(
            \frac{\rho^{\ast}b(V_d)\cap[M_{0,1}(d,X)]}{
                    e_{{\tinyBbb C}^{\times}}(F_0/M_d)}
                     \right)\;
      =\; \frac{ (i_{X_0}^{\ast}\varphi_{\ast}b({\cal U}_d))\,
                  \cap [X_0] }{ e_{{\tinyBbb C}^{\times}}(X_0/W_d) }\,,
           \hspace{1ex}\mbox{denoted}\hspace{1ex}
            \frac{\Theta_d}{e_{{\tinyBbb C}^{\times}}(X_0/W_d)}\,,
              \\[2ex]
   & & & =
     &  g_{\ast}\left(
           \sum_s\,
             \frac{ (\,
               i_{E_{0s}}^{\ast}\,
                 g^{\ast}\,i_{X_0}^{\ast}\varphi_{\ast}b({\cal U}_d)\,
                    )\,
                \cap\,[E_{0s}] }{
                 e_{{\tinyBbb C}^{\times}}(E_{0s}/Quot_{(d)}) }
                  \right)\,,
         \hspace{1ex}\mbox{denoted}\hspace{1ex}
         g_{\ast}\left(
           \sum_s\,
           \frac{\Xi_{d,s}}{
             e_{{\tinyBbb C}^{\times}}(E_{0s}/Quot_{(d)})}
                  \right)\,,
 \end{array}
$$
where $\alpha=c_1({\cal O}_{\CPscriptsize^{\infty}})(1)$ is
 the generator for $H^{\ast}_{{\scriptsizeBbb C}^{\times}}(\pt)$.
On the other hand, one has the intersection numbers and their generating function
$$
 \begin{array}{lllll}
  K_d   & =  & K_d^{V,\,b}   & =
    & \int_{M_{0,0}(d,X)}\,b(V_d)\,,  \\[1.2ex]
  \Phi  & =  & \Phi^{V,\,b}   &  =  &  \sum_d\,K_d\,e^{d\cdot t}\,.
 \end{array}
$$
In the good cases, $K_d$ and $\Phi$ can be obtained from
 $A_d$ and $A(t)$ by appropriate integrals of the form
 $\int_X\,e^{-H\cdot t/\alpha}\,A_d$,
 where $H$ is the hyperplane class on $Y$ restricted to $X$,
 e.g.\ [L-L-Y1, III$\,$: Theorem 3.12].
This integral can be turned into an integral on $E_0\,$:

{\small
$$
 \int_X\,\tau^{\ast}\,e^{H\cdot t}\cap A_d\;
  =\; \int_{Y_0}\,
        e^{H\cdot t}\,
         \cap\,
        g_{\ast}\left(
         \sum_s\,
          \frac{\Xi_{d,s}}{
              e_{{\tinyBbb C}^{\times}}(E_{0s}/Quot_{(d)})}
                \right)\;
  =\; \sum_s\,
        \int_{E_{0s}}\,
          \frac{g^{\ast}e^{H\cdot t}\,\cap\,\widehat{\Xi}_{d,s}}{ 
                 e_{{\tinyBbb C}^{\times}}(E_{0s}/Quot_{(d)})}\,,
 $$
{\normalsize where}} 
$\widehat{\Xi}_{d,s}$ is the Poincar\'{e} dual of $\Xi_{d,s}$
with respect to $[E_{0,s}]$.
As will be discussed in Sec.\ 3.1, $E_{0s}$ is a flag manifold
fibred over $X$ and, hence, $g^{\ast}e^{H\cdot t}$ can be read
off from the natural fibration of flag manifolds
$E_{0s}\rightarrow X$.

Following [L-L-Y1, III$\,$: Sec.\ 5.4], in the case that $b=1$
the above integral is reduced to the integral
$$
 \sum_s\,\int_{E_{0s}}\,
  \frac{g^{\ast}\psi^{\ast}e^{\kappa\cdot\zeta}}{
              e_{{\tinyBbb C}^{\times}}(E_{0s}/\Quot_{(d)})}\,,
$$
where $\kappa$ is the hyperplane class in $W_d$.
In this article, we work out all the equivariant Euler classes
 $e_{{\tinyBbb C}^{\times}}(E_{0s}/\Quot_{(d)})$ 
and hence this integral.

\bigskip

\noindent $\bullet$
{\bf Conventions and notation.}
\begin{itemize}
 \item [(1)]
  All the dimensions are {\it complex} dimensions unless otherwise
  noted.

 \item [(2)]
  The $S^1$-actions involved in this article are induced from
  ${\Bbb C}^{\times}$-actions and both have the same fixed-point
  locus. In many places, it is more convenient to phrase things
  in term of ${\Bbb C}^{\times}$-action and we will not distinguish
  the two actions when this ambiguity causes no harm.

 \item [(3)]
  A locally free sheaf and its associated vector bundle are denoted
  the same.

 \item [(4)]
 An $I\times J$  matrix whose $(i,j)$-entry is $a_{ij}$ is denoted
  by $(a_{ij})_{i,j}$ when the position of an entry is
  emphasized and
  by $[a_{ij}]_{I\times J}$ when the size of the matrix is emphasized.

 \item [(5)]
  From Section 2 on, the smooth curve $C$ will be $\CP^1$
  unless other noted.
\end{itemize}

\bigskip

\bigskip

\section{The $S^1$-action on Quot-schemes.}

Let ${\cal E}^n$ be a trivialized trivial bundle of rank $n$ over $C$.
The $S^1$-action on the Quot-schemes $\Quot({\cal E}^n)$,
 the topology of the $S^1$-fixed-point components, and
 the $S^1$-weights of the normal bundle to these components
 are studied in this section.

\bigskip

\subsection{The $S^1$-fixed-point components.}

We recall first two basic facts that will be needed in the discussion. 

\bigskip

\noindent
{\bf Fact 2.1.1 [modules over P.I.D.].} (Cf.\ [Ja].) {\it
 \begin{itemize}
  \item [{\rm (1)}]
   Let $D$ be a principal ideal domain and $D^{\,\oplus k}$
    be a free module of rank $k$ over $D$.
   Then any submodule of $D^{\,\oplus k}$ is free with basis
    of $\le k$ elements.

  \item [{\rm (2)}]
   If $A\in M_{k\times k}(D)$ be an $k\times k$ matrix with
    entries in $D$, then $A$ is equivalent to a diagonal matrix
    $\Diag\{\,d_1,\,\ldots,\,d_s,\, 0,\,\ldots,\,0\,\}$
    for some $s$, where $d_i\ne 0$ and $d_i|d_j$ if $i\le j$.
    {\rm (}Recall that $A_1,\,A_2\in M_{k\times k}(D)$ are called
     equivalent if $A_2=PA_1Q$ for some invertible
     $P,\,Q\in M_{k\times k}(D)$.{\rm )}
 \end{itemize}
} 

\bigskip

Recall the embedding $S^1=U(1)\hookrightarrow{\Bbb C}^{\times}$,
 which acts on $C={\Bbb C}\cup\{\infty\}$ via $z\rightarrow tz$.
 This lifts to an $S^1$-action (i.e.\ a linearization) on
 the trivialized trivial bundle
 ${\cal E}^n\simeq {\cal O}_C\otimes{\Bbb C}^n$
 given by $(z, v)\mapsto (t\cdot z, v)$.
This induces then an $S^1$-action ${\cal S}\mapsto t\cdot{\cal S}$
 on the set of coherent subsheaves ${\cal S}$ in ${\cal E}$ by
 pulling back local sections$\,$:
 $(t\cdot s)(z)=s(t^{-1}z)$, where $s\in {\cal S}(U)$ and
 $t\cdot s\in (t\cdot{\cal S})(t\cdot U)$
 with $U$ an open subset in $C$.
Since each subsheaf in ${\cal E}^n$ corresponds to a point
 in the Quot-scheme $\Quot({\cal E}^n)$, this gives an $S^1$-action
 $\Quot({\cal E}^n)$.
(Cf.\ [Ak], [Ch], and [Str].)

When restricted to the set of rank-$r$ subbundles in ${\cal E}^n$,
 each holomorphic subbundle in ${\cal E}^n$ corresponds
 to a holomorphic map $f$ from $C$ to a Grassmannian manifold
 $\Gr_r({\Bbb C}^n)$ and
the above $S^1$-action is the $S^1$-action on
 $\Hom(C,\Gr_r({\Bbb C}^n))$
 given by $f\mapsto t\cdot f := f\circ t^{-1}$,
cf.\ [Ak].

In the following, we first characterize the $S^1$-fixed-point in
 $\Quot({\cal E}^n)$ and then give a description of the topology
 of the $S^1$-fixed-point components in $\Quot({\cal E}^n)$.

\bigskip

\noindent
{\bf Lemma 2.1.2 [coherent subsheaf of ${\cal E}^n$].} {\it
 Any coherent subsheaf ${\cal V}$ of ${\cal E}^n$ is locally free.
} 

\bigskip

\noindent {\it Proof.}
 Since any torsion section of ${\cal V}$ is supported on a divisor,
  that support must be contained in an affine chart of the form
  $C-\{\pt\}=\Spec{\Bbb C}[u]$.
 Since ${\cal E}^n$ is globally trivial, ${\cal E}^n|_U$ is
  the sheaf associated to a free ${\Bbb C}[u]$-module $M_U$
  of rank $n$.
 Thus, ${\cal V}|_U$ is the sheaf associated to a submodule
  $M_U^{\prime}$ of $M_U$. Since ${\Bbb C}[u]$ is a principal ideal
  domain, $M_U^{\prime}$ must be free also. This shows that
  ${\cal V}|_U=(M_U^{\prime})^{\sim}$, and hence ${\cal V}$, is
  torsion-free.
 Since a torsion-free coherent sheaf on a smooth curve must be
  locally free, this concludes the lemma.

\noindent\hspace{14cm}$\Box$

\bigskip

\noindent
{\bf Lemma 2.1.3
 [$S^1$-fixed-point $=$ ${\Bbb C}^{\times}$-fixed-point].} {\it
 A coherent subsheaf ${\cal V}$ of ${\cal E}^n$ on $C$
 is $S^1$-invariant if and only if it is ${\Bbb C}^{\times}$-invariant.
} 

\bigskip

\noindent
{\it Proof.}
Only need to show the only-if part.
Let ${\cal V}$ be a rank-$r$ $S^1$-invariant subsheaf in ${\cal E}^n$.
Then ${\cal V}$ is locally free from Lemma 2.1.2 and hence
 there exists an $S^1$-invariant open dense subset
 $U\subset C-\{0,\infty\}$ such that ${\cal V}|_U$ is realized
 as a holomorphic rank-$r$ subbundle of ${\cal E}^n|_U$ and hence
 as a holomorphic map from $U$ into a Grassmannian manifold
 $\Gr_r({\Bbb C}^n)$.
Since ${\cal V}|_U$ is also $S^1$-invarant, this map factors via
 $U\rightarrow U/S^1\rightarrow\Gr_r({\Bbb C}^n)$.
Since $U/S^1$ is a union of open real line segments, holomorphicity
 implies then that any such map must a constant map.
This implies that ${\cal V}|_{C-\{0,\infty\}}$
 is indeed a constant subsheaf in
 ${\cal E}^n|_{C-\{0,\infty\}}$
 with respect to the trivialization of ${\cal E}^n$.
This shows that ${\cal V}$ is also ${\Bbb C}^{\times}$-invariant.

\noindent\hspace{14cm} $\Box$

\bigskip

The following lemma strengthens
 Statement (2) of  Fact 2.1.1   
 in the case of ${\Bbb C}^{\times}$-invariant submodules in
 ${\Bbb C}[z]^{\,\oplus k}$.

\bigskip

\noindent
{\bf Lemma 2.1.4 [${\Bbb C}^{\times}$-invariant submodule].} {\it
 Let $D={\Bbb C}[z]$, $A=A(z)\in \GL(l,{\Bbb C}[z])$
  be an invertible $l\times l$-matrix with entries in ${\Bbb C}[z]$,
  in Fact 2.1.1. 
 If, furthermore, the column vectors of $A(tz)$ generate the same
  ${\Bbb C}[z]$-module for all $t\in {\Bbb C}^{\times}$,
  then there exist invertible $P\in\GL(l,{\Bbb C})$ and
  $Q(z)\in \GL(l,{\Bbb C}[z])$ such that $d_i=z^{\alpha_i}$ in
  Fact 2.1.1 and 
  $$
   A(z)\;
    =\;P\,\Diag\{\,z^{\alpha_1},\,\ldots,\,z^{\alpha_l}\,\}\,Q(z)\,,
  $$
  where $0\le \alpha_1 \le\,\cdots\,\le \alpha_l\,$.  
} 

\bigskip

\noindent
{\it Proof.}
By a sequence of elementary column transformations (e.g.\ [Ja]),
 which correspond to multiplications from the right by a sequence
 of elementary matrices in $\GL(l, {\Bbb C}[z])$, together with
 permutations of rows, which corresponds to a multiplication from
 the left by a sequence of matrices in $\GL(l,{\Bbb C})$,
 one can render $A(z)$ into a lower triangular form
 $B(z)=(b_{ij}(z))_{i,j}$ such that 
 \begin{itemize}
  \item [{\rm (1)}]
   $b_{ij}(z)=0$ for all $i<j$,

  \item [{\rm (2)}]
   $\degree b_{ii}(z) \le \degree b_{i+1,i+1}(z)$ for all $i$, and

  \item [{\rm (3)}]
  $\degree b_{ij}(z) < \degree b_{ii}(z)$ for all $i>j$,
 \end{itemize}
 where $\degree(\,\cdot\,)$ is the degree of the polynomial
 $(\,\cdot\,)$ with respect to the variable $z$ and
 $\degree(0)=-\infty$ by convention.

The assumption that the column vectors of $A(tz)$ generate the same
 ${\Bbb C}[z]$-module for all $t\in{\Bbb C}^{\times}$ implies that
 the column vectors of $B(tz)$ generate the same ${\Bbb C}[z]$-module
 as the module generated by the column vectors of $B(z)$
 for all $t\in{\Bbb C}^{\times}$. In terms of matrices, this is
 equivalent to the existence of $\widehat{Q}(z,t)\in\GL(k,{\Bbb C}[z])$
 such that $B(tz)=B(z)\widehat{Q}(z,t)$, $t\in{\Bbb C}^{\times}$. 
 The fact that both $B(tz)$ and $B(z)$ are lower triangular implies
 that $\widehat{Q}(z,t)$ is also lower triangular.

On the other hand, $\degree b_{ij}(tz)=\degree b_{ij}(z)$
 for all $i,j$. This puts a strong constraint in the form of
 $B(z)$ in order that $B(tz)=B(z)\widehat{Q}(z,t)$ always admits
 a solution for $\widehat{Q}(z,t)$ in $\GL(l,{\Bbb C}[z])$.
Together with the Inequality (3) above:
 $\degree b_{ij}(tz)< \degree b_{ii}(z)$ for all $i>j$, and 
 a tedious yet straightforward induction argument, one can shows
 that $B(z)$ must be of the form
 $$
  B(z)\; =\; B(1)\,\Diag\{\,z^{\alpha_1},\,\ldots,\,z^{\alpha_l}\,\}
 $$
 with $0\le \alpha_1\le\,\ldots\,\le \alpha_l$ and
 $B(1)_{ij}=0$ if $i<j$ or $\alpha_i=\alpha_j$.

This proves the lemma.

\noindent\hspace{14cm}$\Box$

\bigskip

\noindent
{\bf Proposition 2.1.5 [$S^1$-fixed coherent subsheaf].} {\it
 Let ${\cal V}$ be a rank $r$ coherent subsheaf of ${\cal E}^n$
  on $C$.
  Then ${\cal V}$ is a locally free ${\cal O}_C$-module.
 When ${\cal V}$ is in addition $S^1$-invariant, then
  ${\cal V}$ determines a unique enlarged sheaf $\widehat{\cal V}$
  such that
 \begin{itemize}
  \item [{\rm (1)}]
   $\widehat{\cal V}$ is a constant subsheaf in the globally
   trivialized ${\cal E}^n$ of the same rank $r$ as ${\cal V}$,
   {\rm (}thus
     $\widehat{\cal V}\simeq {\cal O}_C^{\,\oplus r}\,${\rm )}.

  \item [{\rm (2)}]
   ${\cal V}$ is a subsheaf of $\widehat{\cal V}$.

  \item [{\rm (3)}]
   Let
    $\{\,U_0=C-\{\infty\}=\Spec{\Bbb C}[z]\,,\,
         U_{\infty}=C-\{0\}=\Spec{\Bbb C}[w]\,\}$
    be an atlas of affine charts on $C$ .
   Then
    there exists a constant re-trivialization
     $$
      \widehat{\cal V}|_{U_0}\;
       =\; {\cal O}|_{U_0}^{\hspace{1ex}\oplus r}\;
       =\; ({\Bbb C}[z]^{\,\oplus r})^{\,\sim}
     $$
    such that
     $$
      {\cal V}|_{U_0}\;
       =\; ({\Bbb C}[z]\,z^{\alpha_1}\,\oplus\,\cdots\,
            \oplus\, {\Bbb C}[z]\,z^{\alpha_r})^{\,\sim}
      \hspace{1em}\mbox{with}\hspace{1em}
      0 \le \alpha_1 \le\,\cdots\,\le \alpha_r
     $$
    with respect to this new trivialization,
    where $(\,\cdot\,)^{\,\sim}$ is the sheaf of modules over
    the affine scheme $U=\Spec R$ in question associated to
    the $R$-module $(\,\cdot\,)$, cf.\ {\rm [Ha]}.
   Similarly for $\widehat{\cal V}|_{U_{\infty}}$ and
    ${\cal V}|_{U_{\infty}}$.
    {\rm (}Corresponding to
            $0\le\beta_1\le\,\ldots,\,\le\beta_r$.{\rm )}
\end{itemize}
}  

\bigskip

\noindent
{\it Remark 2.1.6.} 
 In other words, the local diagonal form of ${\cal V}$ on
  an affine chart can be made compatible with the fixed
  trivialization of ${\cal E}^n$.
 The sheaf ${\cal V}$ can be thought of as obtained by gluing
  the two indenpendent pieces, ${\cal V}|_{U_0}$ and
  ${\cal V}|_{U_{\infty}}$, on affine charts via
  an isomorphism
  $$
   ({\cal V}|_{U_0})|_{U_0\cap U_{\infty}}\,
    \simeq\, ({\Bbb C}[z,z^{-1}]^{\,\oplus r})^{\,\sim}\,
    \stackrel{\hspace{2ex}z\leftrightarrow w^{-1}}{\simeq}\,
       ({\Bbb C}[w^{-1},w]^{\,\oplus r})^{\,\sim}\,
    \simeq\, ({\cal V}|_{U_{\infty}})|_{U_0\cap U_{\infty}}\,.
  $$

\bigskip

\noindent
{\it Proof of Proposition 2.1.5.}
 For Claim (1) and Claim (2).
  Since ${\cal V}$ is an $S^1$-fixed subsheaf in ${\cal E}^n$,
   ${\cal V}|_{C-\{0,\infty\}}$ admits a unique trivial extension
   to a subsheaf of ${\cal E}^n$ on the whole $C$.
   By construction, it has the same rank as ${\cal V}$.
  We shall choose $\widehat{\cal V}$ to be this extension sheaf
   of ${\cal V}_{C-\{0,\infty\}}$.
  If ${\cal V}$ is not contained in $\widehat{\cal V}$ as
   a subsheaf, then there exists an affine chart $U$ of $C$
   such that ${\cal V}|_U$ has a section $s$ not contained in
   $\widehat{\cal V}$.
  Since $\widehat{\cal V}|_U={\cal V}|_U$, this implies that
   $s$ must restrict to the zero-section when localized to
   $U-\{0,\infty\}$.
   In other words, it is a torsion section.
   This contradicts with Lemma 2.1.2, which says that ${\cal V}$
   must be torsion-free.
  Consequently, ${\cal V}$ must be a subsheaf of $\widehat{\cal V}$
   as well.

 For Claim (3).
  Recall Lemma 2.1.4, with $l$ replaced by $r$.
  Since the right multiplication of $A(z)$ by matrices in
  $\GL(r,{\Bbb C}[z])$ does not change the ${\Bbb C}[z]$-module
  generated by the column vectors of $A(z)$, while the left
  multiplication by a constant matrix in $\GL(r,{\Bbb C})$
  can by interpreted as a change of coordinates without
  changing the notion of being a constant section in the
  associated sheaf, this concludes Claim (3) and
  hence the proposition.

\noindent\hspace{14cm}$\Box$

\bigskip

\noindent
{\it Remark 2.1.7.}
 Note the above proposition says that both
  ${\cal V}|_{U_0}$ and ${\cal V}|_{U_{\infty}}$ admit
  diagonalizations by constant global sections in ${\cal E}^n$,
  but in general these two sets of diagonalizing constant sections
  are different.
 This is all right. Indeed for any two such trivializations,
  one over $U_0$ and the other over $U_{\infty}$, the localizations
  of both to $C-\{0,\infty\}$ are isomorphic to the free
  ${\Bbb C}[z,z^{-1}]$-module of rank $r$ and hence they glue
  together to form an ${\Bbb C}^{\times}$-fixed coherent sheaf
  on $C$.

\bigskip

The remaining problem is to decide when two diagonalized
 forms of ${\cal O}_{U_0}$-modules 
 (resp.\ ${\cal O}_{U_{\infty}}$-modules) of rank $r$
 determine the same submodule in ${\cal V}|_{U_0}$
 (resp.\ ${\cal V}|_{U_{\infty}}$).
To determine this, let two diagonal forms of ${\Bbb C}[z]$-modules
 be given by
 $$
  B_1(z)\; =\; B_1(1)\Diag\{z^{\alpha_1},\,\ldots,\,z^{\alpha_r}\}
  \hspace{1em}\mbox{and}\hspace{1em}
  B_2(z)\; =\; B_2(1)\Diag\{z^{\alpha_1},\,\ldots,\,z^{\alpha_r}\}
 $$
 respectively.
Then $B_1(z)$ and $B_2(z)$ determine the same ${\Bbb C}[z]$-module
 if and only if there exists a $Q(z)\in \GL(r,{\Bbb C}[z])$ such
 that $B_1(z)Q(z)=B_2(z)$.
From this, one obtains that
 $$
 \begin{array}{rcl}
   Q(z) & =
    &  \Diag\{z^{-\alpha_1},\,\ldots,\,z^{-\alpha_r}\}\,
        B_1(1)^{-1}\,B_2(1)\,
        \Diag\{z^{\alpha_1},\,\ldots,\,z^{\alpha_r}\} \\[.6ex]
    & = & \Diag\{z^{-\alpha_1},\,\ldots,\,z^{-\alpha_r}\}\,
           B\,
          \Diag\{z^{\alpha_1},\,\ldots,\,z^{\alpha_r}\} \\[.6ex]
    & = & (\, z^{-\alpha_i+\alpha_j}\,b_{ij} \,)_{i,j}\;
          \in \; \GL(r, {\Bbb C}[z])\,,
 \end{array}
 $$
 where $B=B_1(1)^{-1}B_2(1)=(\,b_{ij}\,)_{i,j}$.
 This implies that $b_{ij}=0$ if $\alpha_i>\alpha_j$.
Consequently, $Q(z)$ is a {\it block upper triangular matrix, 
 whose block form is determined by the multiplicity of elements
 in $(\alpha_1,\,\ldots,\,\alpha_r)$. }
 (For example, if this sequence is $(1, 1, 4, 4, 4, 7)$, then
  the corresponding block upper triangular matrix will have in
  the diagonal $2\times 2$-, $3\times 3$-, and $1\times 1$-blocks.)
Rephrased in a more geometric way,
 {\it $B_1(z)$ and $B_2(z)$ determine the same submodule
      if and only if they correspond to the same flag.}
Explicitly, the flag associated to 
 $B(z)=B(1)\Diag\{z^{\alpha_1},\,\ldots,\, z^{\alpha_r}\}$
 is given as follows.

Let $B(1)=(u_1,\,\ldots,u_r)$ be the column vectors of $B(1)$
 and suppose that
 $$
  \alpha_1=\cdots=\alpha_{j_1}<
   \alpha_{j_1+1}=\cdots=\alpha_{j_2} < \cdots
    < \alpha_{j_s+1}=\cdots=\alpha_r\,,
 $$
 then the flag associated to $B(z)$ is given by
{\small
 $$ \hspace{-1em}
  (\,\Span\{u_1,\,\ldots,\,u_{j_1}\}\;
      \subset\;\Span\{u_1,\,\ldots,\,u_{j_2}\}\;
      \subset\; \cdots\;
      \subset\;\Span\{u_1,\,\ldots,\,u_{j_s}\}\;
      \subset\; {\smallBbb C}^r\,)
     \;\in\; \Fl_{j_1,\,\ldots,\,j_s}({\smallBbb C}^r)\,,
 $$
 {\rm where}} 
 note that the last ${\Bbb C}^r$ should be identified with
 $\widehat{\cal V}\in \Gr_r({\Bbb C}^n)$.

\bigskip

\noindent
{\bf Definition 2.1.8 [admissible pair of sequences].} {
 Recall the Hilbert polynomial $P=P(t)=(n-r)t+d+(n-r)$.
 Then $(\alpha_1,\,\ldots,\,\alpha_r\;;\;\beta_1,\,\ldots,\,\beta_r)$
  is called an {\it admissible pair of sequences} with respect
  to $P(t)$ if it satisfies 
 \begin{itemize}
  \item [(1)]
   $0\le \alpha_1\le\,\ldots\,\le \alpha_r$, \hspace{1ex}
   $0\le\beta_1\le\,\ldots\,\le \beta_r$, and 
   
  \item [(2)]
   $(\alpha_1+\,\ldots\,+\alpha_r)+(\beta_1+\,\ldots\,+\beta_r)=d\,$.
 \end{itemize}
} 

\bigskip

From the above discussions and the fact that, for any element
 in $\Fl_{j_1,\,\ldots,\,j_s}({\Bbb C}^r)$, one can always
 construct a $B(z)$ in the above form such that $B(1)$
 is mapped to that flag by the above correspondence,
 one concludes the following proposition.

\bigskip

\noindent 
{\bf Theorem 2.1.9 [topology of $S^1$-fixed-point locus].} {\it
  Let
  $(\alpha_1,\,\ldots,\,\alpha_r\;;\;\beta_1,\,\ldots,\,\beta_r)$
   be an admissible pair of sequences of non-negative integers,
  $\Fl_{j_1,\ldots,\,j_s, r}({\Bbb C}^{n})$ and
   $\Fl_{j^{\prime}_1,\ldots,\,j^{\prime}_{s^{\prime}}, r}(
                                                  {\Bbb C}^{n})$
   be the flag manifold associated to the multiplicity of elements 
   in $(\alpha_1,\,\ldots,\,\alpha_r)$ and
   $(\beta_1,\,\ldots,\,\beta_r)$ respectively,
   as discussed above.
  Let
  $$
   \Fl_{j_1,\ldots,\,j_s, r}({\Bbb C}^{n})
                         \rightarrow \Gr_r({\Bbb C}^n)
    \hspace{1em}\mbox{and}\hspace{1em}
   \Fl_{j^{\prime}_1,\ldots,\,j^{\prime}_{s^{\prime}}, r}(
                                                  {\Bbb C}^{n})
                          \rightarrow \Gr_r({\Bbb C}^n)
  $$
  be the natural projections.
  Then the subset
  $F_{\alpha_1,\,\ldots,\,\alpha_r\,;\,\beta_1,\,\ldots,\,\beta_r}$
  of the $S^1$-fixed-point locus that is associated to
  $(\alpha_1,\,\ldots,\,\alpha_r\,;\,\beta_1,\,\ldots,\,\beta_r)$
  is connected and is given by the fiber product
  $$
   \Fl_{j_1,\ldots,\,j_s, r}({\Bbb C}^{n})\;
     \times_{Gr_r({\scriptsizeBbb C}^n)}\;
   \Fl_{j^{\prime}_1,\ldots,\,j^{\prime}_{s^{\prime}}, r}(
                                                  {\Bbb C}^{n})\,.
  $$
} 

\bigskip

\noindent 
{\it Remark 2.1.10.}
 \begin{itemize}
  \item [{\rm (1)}]
   The base $\Gr_r({\Bbb C}^n)$ corresponds to the choices
    of $\widehat{\cal V}$.
   The fiber over a point in the base is the product of
    two flag manifolds that gives all possible $S^1$-fixed subsheaves
    ${\cal V}$ of $\widehat{\cal V}$ that have the specified
    Hilbert polynomial of ${\cal E}^n\!/{\cal V}$ associated to
    $(\alpha_1,\,\ldots,\,\alpha_r\,;\,\beta_1,\,\ldots,\,\beta_r)$.

  \item [{\rm (2)}]
   The dimension of components of $S^1$-fixed-point locus varies
    from component to component. When the $S^1$-fixed-point locus
    is non-empty, the dimension of each component is bounded below
    by the dimension of the target Grassmannian manifold
    $\Gr_r({\Bbb C}^n)$ that one starts with.
   The above expression implies that the only case that a fixed-point
    component has the dimension the same as that
    of $\Gr_r({\Bbb C}^n)$ is when that component itself is
    homeomorphic to $\Gr_r({\Bbb C}^n)$. This happens exactly
    when $\alpha_1=\cdots=\alpha_r$ and $\beta_1=\cdots=\beta_r$.
    Such $(\alpha_1,\,\cdots,\alpha_r\;;\;\beta_1,\cdots\,,\beta_r)$
    is admissible only for special Hilbert polynomials.
\end{itemize}

\bigskip

\subsection
{\bf The $S^1$-weight system of the tangent space of Quot-scheme
     at an $S^1$-fixed-point component. }

After recalling some related facts in the preparatory remarks,
we compute the $S^1$-weights and their multiplicities of
the tangent space of the Quot-scheme $\Quot({\cal E}^n)$
at an $S^1$-fixed-point.

\bigskip

\begin{flushleft}
{\bf Preparatory remarks.}
\end{flushleft}
Recall (cf.\ [Ch], [H-L], and [Kol]) that the tangent space of
 Quot-scheme at a point is given by 
 $$
  T_{({\cal E}^n\rightarrow
   {\cal E}^n\!/\mbox{\raisebox{-.4ex}{\scriptsize${\cal V}$}}
   )}\Quot_P({\cal E}^n)\;
  \simeq\; \Hom_{{\cal O}_C}(\,
            {\cal V}\,,\,
            {\cal E}^n\!/\mbox{\raisebox{-.4ex}{${\cal V}$}}\,)\,.
 $$
When
 $({\cal E}^n \rightarrow
       {\cal E}^n\!/\mbox{\raisebox{-.4ex}{${\cal V}$}})$
 is an $S^1$-fixed-point, $S^1$ acts both on ${\cal V}$ and
 ${\cal E}^n\!/\mbox{\raisebox{-.4ex}{${\cal V}$}}$.
 The $S^1$-action on
 $T_{({\cal E}^n\rightarrow
    {\cal E}^n\!/\mbox{\raisebox{-.4ex}{\scriptsize${\cal V}$}}
        )}\Quot_P({\cal E}^n)$
 is translated to the $S^1$-action on
 $\Hom_{{\cal O}_C}(\,
     {\cal V}\,,\,
      {\cal E}^n\!/\mbox{\raisebox{-.4ex}{${\cal V}$}}\,)$
 by conjugations: $f\,\mapsto \, t\cdot f\cdot t^{-1}$ for $t\in S^1$.

Recall the inclusion of $S^1$-invariant subsheaves
 ${\cal V}\subset \widehat{\cal V}$ in ${\cal E}^n$.
 One thus has a natural morphism
 ${\cal E}^n\!/\mbox{\raisebox{-.4ex}{${\cal V}$}}
  \rightarrow
  {\cal E}^n\!/\mbox{\raisebox{-.4ex}{$\widehat{\cal V}$}}$.
 Since $\widehat{\cal V}$ is a constant rank-$r$ subbundle in
  ${\cal E}^n$,
  ${\cal E}^n\!/\mbox{\raisebox{-.4ex}{$\widehat{\cal V}$}}$
  is a rank-$(n-r)$ trivial bundle on $C$. 
Since
 $({\cal E}^n\!/\mbox{\raisebox{-.4ex}{${\cal V}$}})|
                                             _{C-\{0,\infty\}}
   \stackrel{\sim}{\rightarrow}
   ({\cal E}^n\!/\mbox{\raisebox{-.4ex}{$\widehat{\cal V}$}})|
                                             _{C-\{0,\infty\}}$
  from the restriction of the above morphism and
 the restriction to the stalks
  $$
   ({\cal E}^n\!/\mbox{\raisebox{-.4ex}{${\cal V}$}})_0
    \rightarrow
   ({\cal E}^n\!/\mbox{\raisebox{-.4ex}{$\widehat{\cal V}$}})_0
   \hspace{2em}
   (\;\mbox{resp.}\hspace{1ex}
   ({\cal E}^n\!/\mbox{\raisebox{-.4ex}{$\widehat{\cal V}$}})
                                                      _{\infty}
     \rightarrow
     ({\cal E}^n\!/\mbox{\raisebox{-.4ex}{$\widehat{\cal V}$}})
                                                      _{\infty})
  $$
  at $0$ (resp.\ $\infty$) given by
  {\small
  $$
   \left(\,{\smallBbb C}[z]^{\oplus n}/
     ({\smallBbb C}[z]z^{\alpha_1}\oplus\,
        \cdots\,\oplus{\smallBbb C}[z]z^{\alpha_r}\,
    )\right)
    \otimes_{{\cal O}_C(U_0)}{\cal O}_{C,\,0}\;
   \longrightarrow\;
   {\smallBbb C}[z]^{\oplus (n-r)}  
            \otimes_{{\cal O}_C(U_0)}{\cal O}_{C,\,0}\;
  $$
  {\normalsize (resp.} 
  $$
   \left(\,{\smallBbb C}[z]^{\oplus n}/
     ({\smallBbb C}[z]z^{\beta_1}\oplus\,
        \cdots\,\oplus{\smallBbb C}[z]z^{\beta_r}\,
    )\right)
    \otimes_{{\cal O}_C(U_{\infty})}{\cal O}_{C,\,\infty}\;
   \longrightarrow\;
   {\smallBbb C}[z]^{\oplus (n-r)}  
    \otimes_{{\cal O}_C(U_{\infty})}{\cal O}_{C,\,\infty}\;)
  $$
  {\normalsize are}} 
  surjective, the morphism
  ${\cal E}^n\!/\mbox{\raisebox{-.4ex}{${\cal V}$}}
    \rightarrow
    {\cal E}^n\!/\mbox{\raisebox{-.4ex}{$\widehat{\cal V}$}}$
  is surjective and one has the following split exact sequence
  of torsion-part/locally-free-part decomposition
  $$
   0\;\longrightarrow\;
    \widehat{\cal V}/\mbox{\raisebox{-.4ex}{${\cal V}$}}\;
     \longrightarrow\;
    {\cal E}^n\!/\mbox{\raisebox{-.4ex}{${\cal V}$}}\;
      \longrightarrow\;
    {\cal E}^n\!/\mbox{\raisebox{-.4ex}{$\widehat{\cal V}$}}\;
      \longrightarrow\; 0\,.
 $$
 Since any constant rank-$(n-r)$ subbundle in ${\cal E}^n$
 that is transverse to $\widehat{\cal V}$ is $S^1$-invariant
 and is mapped isomorphically to
 ${\cal E}^n\!/\mbox{\raisebox{-.4ex}{$\widehat{\cal V}$}}$,
 the above decomposition is also $S^1$-equivariant.

\bigskip

\begin{flushleft}
{\bf The $S^1$-action on
  $\Hom_{{\cal O}_C}(\,
    {\cal V}\,,\,
    {\cal E}^n\!/\mbox{\raisebox{-.4ex}{${\cal V}$}}\,)$
  when
  $({\cal E}^n\rightarrow {\cal E}^n\!
        /\mbox{\raisebox{-.4ex}{${\cal V}$}})$
  is an $S^1$-fixed-point.}
\end{flushleft}
The above discussion gives an $S^1$-invariant decomposition of
the tangent space to the Quot-scheme at an $S^1$-fixed-point $\,$:
\begin{eqnarray*}
 \lefteqn{
  \Hom_{{\cal O}_C}(\,{\cal V}\,,\,
       {\cal E}^n\!/\mbox{\raisebox{-.4ex}{${\cal V}$}}\,)\;
  =\; \Hom_{{\cal O}_C}(\,{\cal V}\,,\,
      \widehat{\cal V}/\mbox{\raisebox{-.4ex}{${\cal V}$}}
       \oplus
    {\cal E}^n\!/\mbox{\raisebox{-.4ex}{$\widehat{\cal V}$}}\,)\;}
                                                        \\[.6ex]
  && =\;  \Hom_{{\cal O}_C}(\,{\cal V}\,,\,{\cal F}_0\,)\;
     \oplus\;
     \Hom_{{\cal O}_C}(\,{\cal V}\,,\,{\cal F}_{\infty}\,)\;
     \oplus\;
     \Hom_{{\cal O}_C}(\,{\cal V}\,,\,
     {\cal E}^n\!/\mbox{\raisebox{-.4ex}{$\widehat{\cal V}$}}\,)\,,
\end{eqnarray*}
where ${\cal F}_0$ (resp.\ ${\cal F}_{\infty}$) is the torsion
subsheaf of ${\cal E}^n\!/\mbox{\raisebox{-.4ex}{${\cal V}$}}$
supported at $0$ (resp.\ $\infty$).
We shall now study the three summands in the decomposition and
 their $S^1$-weight system, denoted by $\Wt_{\,1}$, $\Wt_{\,2}$,
 and $\Wt_{\,3}$
 respectively.
Due to the tediousness of the discussion, we itemize the argument
 below.

\bigskip

\begin{flushleft}
{{\it $\bullet$
  The summands
   $\Hom_{{\cal O}_C}(\,{\cal V}\,,\,{\cal F}_0\,)$
   and
   $\Hom_{{\cal O}_C}(\,{\cal V}\,,\,{\cal F}_{\infty}\,)$}$\,$:}
\end{flushleft}
\begin{itemize}
 \item [(1)]
  These two components can be calculated via the restriction
   of the former to $U_0$ and the latter to $U_{\infty}$.
  The problem is reduced then to the study of the group of
   homomorphisms of ${\Bbb C}[z]$-modules and the $S^1$-action
   on it.
  Explicitly,
  $$
   \Hom_{{\cal O}_C}(\,{\cal V}\,,\,{\cal F}_0\,)\;
   =\;\Hom_{{\scriptsizeBbb C}[z]}
       (\,\oplus_{j=1}^r\,{\Bbb C}[z]\cdot z^{\alpha_j}\,,\,
     \oplus_{i=1}^r\,{\Bbb C}[z]\cdot\overline{e}_{0\,i}\,)\,,
  $$
  where $z^{\alpha_i}\cdot\overline{e}_{0\,i}=0$
  for $i=1,\,\ldots,\,r\,$, and
  $$
   \Hom_{{\cal O}_C}(\,{\cal V}\,,\,{\cal F}_{\infty}\,)\;
   =\;\Hom_{{\scriptsizeBbb C}[w]}
     (\,\oplus_{j=1}^r\,{\Bbb C}[w]\cdot w^{\beta_j}\,,\,
     \oplus_{i=1}^r\,{\Bbb C}[w]\cdot\overline{e}_{\infty\,i}\,)\,,
  $$
  where $w^{\beta_i}\cdot\overline{e}_{\infty\,i}=0$
  for $i=1,\,\ldots,\,r$.

 \medskip
 \item [(2)]
 {\it Computation of the weight systems $\Wt_{\,1}$ and $\Wt_{\,2}$}$\,$:

 \medskip
 \item []
 \hspace{-1em}(2.1)\hspace{1ex}
  Realize an element in
  $\oplus_{i=1}^r\,{\Bbb C}[z]\cdot\overline{e}_{0\,i}$
  as a column vector and let
  $$
   f(z)\;=\; (f_{ij}(z))_{i,j}\;
    \in\; \Hom_{{\cal O}_C}(\,{\cal V}\,,\,{\cal F}_0\,)
  $$
  with respect to the local bases
  $(z^{\alpha_1},\,\ldots,\, z^{\alpha^r})$
  and $(\overline{e}_{01},\,\ldots,\,\overline{e}_{0r})$
  for ${\cal V}$ and ${\cal F}_0$ respectively.
  Then $\degree f_{ij}(z)<\alpha_i$ and 
  (cf.\ the similar computation for the weight system $\Wt_{\,3}$ below),
  $$
   (t\cdot f)(z)\;=\; (\, t^{\alpha_j}\,f_{ij}(t^{-1}z)\,)_{i,j}\,,
    \hspace{1em}t\in S^1\,.
  $$
  Thus, the rank-$1$ $S^1$-eigen-spaces in
  $\Hom_{{\cal O}_C}(\,{\cal V}\,,\,{\cal F}_0\,)$
  can be chosen to be generated by 
  $$
   \hspace{-3em}
   \varepsilon_{ij}[\mu^0_{ij}]\;
    = \; \varepsilon_{ij}[\mu^0_{ij}](z)\;
    := \; (f_{kl}(z))_{k,l}\,,
    \hspace{1em}\mbox{where}\hspace{1em}
    f_{kl}(z)\,\;
     =\; \left\{
          \begin{array}{ll}
           0  & \mbox{if $(k,l)\ne (i,j)$,} \\[.6ex]
           z^{\alpha_j-\mu^0_{ij}}
              & \mbox{if $(k,l)=(i,j)$.}
          \end{array}
         \right.
  $$
  whose $S^1$-weight is $\mu^0_{ij}$ that satisfies
  $$
    \alpha_j-\alpha_i\;<\;\mu^0_{ij}\;\le\;\alpha_j\,.
  $$
  From this, one has
  $$
   \Wt_{\,1}\;=\;\bigsqcup_{\,i,\,j=1}^{\,r}\,
     (\,(\,\alpha_j-\alpha_i\,,\,\alpha_j\,]\,\cap\,{\Bbb Z}\,)
  $$
  with the multiplicity of a given integer in the set
  being the number of times it appears in the disjoint union.

 \medskip
 \item []
 \hspace{-1em}(2.2)\hspace{1ex}
  Rewrite $(\alpha_1,\,\ldots,\,\alpha_r\,)$ as
   $$
    \begin{array}{ccccccc}
     0 & \le  & a_1\,(=\alpha_1)
                    & <  & \cdots & < & a_k\,(=\alpha_r)\,,
                                                       \\[.6ex]
       &      & m_1 &    & \cdots &   & m_k
     \end{array}
   $$
   with the multiplicity indicated.
  For an interval $I\subset{\Bbb R}$, let $\chi_I$ be
   the characteristic function
   $\chi_I(x)=1$, if $x\in I$, and $=0$, otherwise.
  Let $\chi^A\,=\,\sum_{i=1}^k\,m_i\,\chi_{(-a_i\,,\,0]}$
   and define $\chi^A_m$ by
   $\chi^A_m(x)\,=\,\sum_{j=1}^k\,m_j\,\chi^A(x-a_j)$.
  Then the multiplicity for $\mu\in \Wt_{\,1}$ is given by
   $\chi^A_m(\mu)$.
     
 \medskip
 \item []
 \hspace{-1em}(2.3)\hspace{1ex}
  Realize an element in
  $\oplus_{i=1}^r\,{\Bbb C}[w]\cdot\overline{e}_{\infty\,i}$
  as a column vector and let
  $$
   g(w)\;=\; (g_{ij}(w))_{i,j}\;
    \in\; \Hom_{{\cal O}_C}(\,{\cal V}\,,\,{\cal F}_{\infty}\,)\,.
  $$
  Then $\degree g_{ij}(w)<\beta_i$ and
  (cf.\ the similar computation for the weight system $\Wt_{\,3}$ below),
  $$
   (t\cdot g)(w)\;=\; (\, t^{-\beta_j}\,g_{ij}(tw)\,)_{i,j}\,,
    \hspace{1em}t\in S^1\,.
  $$
  Thus, the rank-$1$ $S^1$-eigen-spaces in
  $\Hom_{{\cal O}_C}(\,{\cal V}\,,\,{\cal F}_{\infty}\,)$
  can be chosen to be generated by
  $$
   \hspace{-3em}
   \varepsilon_{ij}[\mu^{\infty}_{ij}]\;
    = \; \varepsilon_{ij}[\mu^{\infty}_{ij}](z)\;
    := \; (g_{kl}(z))_{k,l}\,,
    \hspace{1em}\mbox{where}\hspace{1em}
    g_{kl}(w)\,\;
     =\; \left\{
          \begin{array}{ll}
           0  & \mbox{if $(k,l)\ne (i,j)$,} \\[.6ex]
           w^{\beta_j+\mu^{\infty}_{ij}}
              & \mbox{if $(k,l)=(i,j)$.}
          \end{array}
         \right.
  $$
  whose $S^1$-weight is $\mu^{\infty}_{ij}$ that satisfies
  $$
   -\beta_j\;\le\;\mu^{\infty}_{ij}\;<\;\beta_i-\beta_j\,.
  $$
  From this, one has
  $$
   \Wt_{\,2}\;=\;\bigsqcup_{\,i,\,j=1}^{\,r}\,
     (\,[\,-\beta_j\,,\,\beta_i-\beta_j\,)\,\cap\,{\Bbb Z}\,)
  $$
  with the same rule of counting multiplicity as for $\Wt_{\,1}$.

 \medskip
 \item []
 \hspace{-1em}(2.4)\hspace{1ex}
  Rewrite $(\beta_1,\,\ldots,\,\beta_r)$ as
   $$
    \begin{array}{ccccccccccccccc}
     0 & \le  & b_1\,(=\beta_1)
              & <  & \cdots & < & b_l\,(=\beta_r) \\[.6ex]
              &      & n_1 &    & \cdots &   & n_l
    \end{array}
   $$
   with the multiplicity indicated.
  Let $\chi^B\,=\,\sum_{i=1}^l\,n_i\,\chi_{[\,0\,,\,\beta_i\,)}$
   and define $\chi^B_m$ by \newline
   $\chi^B_m(x)\,=\,\sum_{j=1}^l\,n_j\,\chi^A(x+b_j)$.
  Then the multiplicity for $\mu\in \Wt_{\,2}$ is given by
   $\chi^B_m(\mu)$.
\end{itemize}

\bigskip

\begin{flushleft}
{{\it $\bullet$
   The summand
    $\Hom_{{\cal O}_C}(\,{\cal V}\,,\,
   {\cal E}^n\!/\mbox{\raisebox{-.4ex}{$\widehat{\cal V}$}}\,)\,$}:}
\end{flushleft}
\begin{itemize}
 \item [(1)]
  Since ${\cal E}^n\!/\mbox{\raisebox{-.4ex}{$\widehat{\cal V}$}}$
   is represented by a rank-$(n-r)$ constant subbundle
   in ${\cal E}^n$ transverse to $\widehat{\cal V}$, it can be
   further decomposed into a direct sum of constant line subbundles
   in ${\cal E}^n$.
   Since all the bundles involved are constant, the decomposition
   of the quotient 
   ${\cal E}^n\!/\mbox{\raisebox{-.4ex}{$\widehat{\cal V}$}}
    = {\cal O}_C^{\oplus (n-r)}$
   is $S^1$-invariant.
  Recall that $S^1$ acts on ${\cal E}^n$ and hence on
   ${\cal E}^n\!/\mbox{\raisebox{-.4ex}{$\widehat{\cal V}$}}$
   via the trivial linearization.
  With respect to this decomposition, one has 
  $$
   \Hom_{{\cal O}_C}(\,{\cal V}\,,\,
     {\cal E}^n\!/
       \mbox{\raisebox{-.4ex}{$\widehat{\cal V}$}}\,)\;
   =\; \Hom_{{\cal O}_C}(\,{\cal V}\,,\,
                             {\cal O}_C^{\,\oplus (n-r)}\,)\;
   =\; H^0(C,{\cal V}^{\vee})^{\oplus (n-r)}\,.                   
  $$

 \medskip
 \item []
  {\it Remark.}
  A connected component of the fixed-point locus can be
   stratified by subsets labelled by the isomorphism classes
   of vector bundles associated to the $S^1$-invariant subsheaves
   ${\cal V}$ in ${\cal E}^n$.
  There can be more than one strata for a connected component. 
\end{itemize}

\begin{itemize}
 \item [(2)]
 {\it Computation of the weight system $\Wt_{\,3}$}$\,$:

 \medskip
 \item []
 \hspace{-1em}(2.1)\hspace{1ex}
  Recall the $S^1$-invariant decompositions
  $$
   \Hom_{{\cal O}_C}(\,{\cal V}\,,\,
     {\cal E}^n\!/
       \mbox{\raisebox{-.4ex}{$\widehat{\cal V}$}}\,)\;
   =\; \Hom_{{\cal O}_C}(\,{\cal V}\,,\,
                             {\cal O}_C^{\,\oplus (n-r)}\,)\;
   =\; \Hom_{{\cal O}_C}(\,{\cal V}\,,\,
                          {\cal O}_C\,)^{\,\oplus (n-r)}\,.
  $$
  The existence of such $S^1$-invariant decomposition implies
  that the sought-for $S^1$-weight system for
  $\Hom_{{\cal O}_C}(\,{\cal V}\,,\,
      {\cal E}^n\!/
            \mbox{\raisebox{-.4ex}{$\widehat{\cal V}$}}\,)$
  consists of $(n-r)$-many copies of the $S^1$-weight system for
  $\Hom_{{\cal O}_C}(\,{\cal V}\,,\,{\cal O}_C\,)$.

 \medskip
 \item []
  \hspace{-1em}(2.2)\hspace{1ex} 
  Let $f\in\Hom_{{\cal O}_C}(\,{\cal V}\,,\,{\cal O}_C\,)$.
  Then, as a morphism of sheaves on $C$, $f$ is described by
  a pair
  {\small
  $$
   \begin{array}{c}
    (f_0, f_{\infty})\;\in\;
       \Hom_{{\cal O}_C(U_0)}
                  (\,{\cal V}(U_0)\,,\,{\cal O}_C(U_0)\,)\,
     \times\,
     \Hom_{{\cal O}_C(U_{\infty})}(\,{\cal V}(U_{\infty})\,,
                              \,{\cal O}_C(U_{\infty})\,) \\[1ex]
    \hspace{3em}
       =\; \Hom_{{\tinyBbb C}[z]}
             (\,\oplus_{i=1}^r
                  {\smallBbb C}[z]\cdot z^{\alpha_i}\,,\,
                                         {\smallBbb C}[z]\,)\,
            \times\,
            \Hom_{{\tinyBbb C}[w]}
              (\,\oplus_{i=1}^r
                   {\smallBbb C}[w]\cdot w^{\beta_i}\,,\,
                                         {\smallBbb C}[w]\,)
   \end{array}
  $$
  {\normalsize such that}} 
  $f_0|_{U_0\cap U_{\infty}}=f_{\infty}|_{U_0\cap U_{\infty}}$.

 \medskip
 \item []
  \hspace{-1em}(2.3)\hspace{1ex}
   Recall the proof of Theorem 5.3 in [Ko], which says in our case
   that {\it the weight system of the tangent bundle at an
   $S^1$-fixed-point depends only on the connected component
   of the fixed-point locus}. 
  Thus, {\it
   to compute the weight system one can choose the $S^1$-invariant
   subsheaf ${\cal V}$ in ${\cal E}^n$ such that
   the two diagonalized local pieces on affine charts $U_0$ and
   $U_{\infty}$ match {\rm (}i.e.\ ${\cal V}$ becomes
   the direct sum of appropriate ideal sheaves in constant line
   subbundles in ${\cal E}^n${\rm )}}.
  From the previous discussions, there are many
   - even continuous families of - such ${\cal V}$.
   However, as will be clear from the explicit expression that 
   the weight system obtained is indeed independent of which
   such ${\cal V}$ is chosen for the computation, as long as
   they belong to the same fixed-point component.
  This gives a consistency check of the method.

 \medskip
 \item []
  \hspace{-1em}(2.4)\hspace{1ex}
  Let ${\cal V}$ be an $S^1$-invariant subsheaf of ${\cal E}^n$
   such that the two local diagonalizations match and 
   suppose that $\alpha_i$ is matched with $\beta_{i^{\prime}}$,
   $i=1,\,\ldots,\,r$.
  Then ${\cal V}$ is decomposed into a direct sum
   $\oplus_{i=1}^r{\cal I}_{\alpha_i,\beta_{i^{\prime}}}$,
   where ${\cal I}_{\alpha_i,\beta_{i^{\prime}}}$ is a subsheaf
   in a constant line subbundle $\simeq{\cal O}_C$ in ${\cal E}^n$
   with the local data as a sheaf of ${\cal O}_C$-module:
   $$
    \begin{array}{cccccccl}
     \mbox{\small on $U_0$}
       & & \mbox{\small on $U_0\cap U_{\infty}$}
       & & \mbox{\small on $U_0\cap U_{\infty}$}
       & & \mbox{\small on $U_{\infty}$}         & \\[.6ex]
     {\Bbb C}[z]\cdot z^{\alpha_i} 
       & \rightarrow & {\Bbb C}[z, z^{-1}] 
       & \stackrel{z\leftrightarrow 1/w}{\longleftrightarrow}
       & {\Bbb C}[w, w^{-1}] 
       & \leftrightarrow    & {\Bbb C}[w]\cdot w^{\beta_{i^{\prime}}}
       & .
    \end{array}
   $$
  The vector bundle associated to ${\cal V}$ is isomorphic
   to $\oplus_{i=1}^r{\cal O}(-\alpha_i-\beta_{i^{\prime}})$
   and $\Hom_{{\cal O}_C}(\,{\cal V}\,,\,{\cal O}_C\,)$
   is further decomposed into an $S^1$-invariant direct sum
   $$
     \Hom_{{\cal O}_C}(\,{\cal V}\,,\,{\cal O}_C\,)\;
     =\; \oplus_{i=1}^r\,
           \Hom_{{\cal O}_C}
            (\,{\cal I}_{\alpha_i, \beta_{i^{\prime}}}\,,\,
                                                  {\cal O}_C\,)\,.
   $$

 \item []
  \hspace{-1em}(2.5)\hspace{1ex}
   For simplicity of notation, we shall drop temporarily
    the indices $i$ and $i^{\prime}$.
   At the level of sheaf morphisms, the data that encodes 
    $f\in\Hom_{{\cal O}_C}
                 (\,{\cal I}_{\alpha,\beta}\,,\,{\cal O}_C\,)$
    is given by a pair
   $$
    \begin{array}{cl}
     (f_0, f_{\infty})\;\in\;
       \Hom_{{\scriptsizeBbb C}[z]}
             ({\Bbb C}[z]\cdot z^{\alpha}\,,\,{\Bbb C}[z])
       \times
       \Hom_{{\scriptsizeBbb C}[w]}
              ({\Bbb C}[w]\cdot w^{\beta}\,,\,{\Bbb C}[w])
                                                 & ,\\[.6ex]
      (\,z^{\alpha}\stackrel{f_0}{\longmapsto} h_0(z)\,,\,
         w^{\beta}\stackrel{f_{\infty}}{\longmapsto} h_{\infty}(w)\,)
    \end{array}
   $$
   such that the following matching condition holds
   $$
    z^{-\alpha}h_0(z)\;=\; w^{-\beta}h_{\infty}(w)
    \hspace{1em}\mbox{under $z\rightarrow 1/w$}\,. 
   $$
  Consequently,
  {\small
   \begin{eqnarray*}
     \lefteqn{
      \Hom_{{\cal O}_C}
        (\,{\cal I}_{\alpha,\beta}\,,\,{\cal O}_C\,)
         \;=\; \{\,(h_0(z)\,,\,h_{\infty}(w))\,|\,
                \degree h_0(z)\le \alpha+\beta
                 \hspace{1ex}\mbox{and}\hspace{1ex}
                   h_{\infty}(w)=w^{\alpha+\beta}h_0(1/w)\,\}
                                                      } \\[.6ex]
    & &  \simeq\; H^0(C,{\cal I}^{\vee})\;
           =\; H^0(C,{\cal O}_C(\alpha+\beta))\,. \hspace{18em}
   \end{eqnarray*}
  } 

 \item []
 \hspace{-1em}(2.6)\hspace{1ex}
  The $S^1$-action on
   $\Hom_{{\cal O}_C}(\,{\cal I}_{\alpha,\beta}\,,\,{\cal O}_C\,)$
   is given by $f\mapsto t\cdot f$, where
   $t\cdot f$ is the composition of the following conjugation of
   $f=(\,h_0(z)\,,\,h_{\infty}(w)\,)\,$:
   {\scriptsize
   $$ \hspace{-13em}
     \begin{array}{lccccccl}
      \mbox{on $U_0\,$:} \\
      & s_0(z)z^{\alpha} & \stackrel{t^{-1}}{\longrightarrow}
        & s_0(tz)(tz)^{\alpha} = t^{\alpha}s_0(tz)\cdot z^{\alpha}
        & \stackrel{f}{\longrightarrow}
        & t^{\alpha}s_0(tz)\,h_0(z)  & \stackrel{t}{\rightarrow}
        & t^{\alpha}s_0(tt^{-1}z)h_0(t^{-1}z)
          = t^{\alpha}s_0(z)h_0(t^{-1}z)              \\[1ex]
     \mbox{on $U_{\infty}\,$:} \\
      & s_{\infty}(w)w^{\beta} & \stackrel{t^{-1}}{\longrightarrow}
      & s_{\infty}(t^{-1}w)(t^{-1}w)^{\beta}
          = t^{-\beta}s_{\infty}(t^{-1}w)\cdot w^{\beta}
        & \stackrel{f}{\longrightarrow}
        & t^{-\beta}s_{\infty}(t^{-1}w)\,h_{\infty}(w)
        & \stackrel{t}{\rightarrow}
        & t^{-\beta}s_{\infty}(t^{-1}tw)h_{\infty}(tw)
          = t^{-\beta}s_{\infty}(w)h_{\infty}(tw)\,.   \\[1ex]
    \end{array}
   $$
   {\normalsize One}} 
   can check directly that if $(f_0,\, f_{\infty})$ satisfies
   the matching condition, then so does
   $(t\cdot f)_0,\, (t\cdot f)_{\infty})$. 
   Consequently,
   $$
    f=(\,h_0(z)\,,\,h_{\infty}(w)\,)\;
     \stackrel{t}{\longrightarrow}\;
     t\cdot f
     = (\,t^{\alpha}\,h_0(t^{-1}z)\,,\,
                      t^{-\beta}\,h_{\infty}(tw)\,)
      \hspace{1ex}\mbox{on}\hspace{1ex}
     \Hom_{{\cal O}_C}
               (\,{\cal I}_{\alpha,\beta}\,,\,{\cal O}_C\,)\,.
   $$
  If $f$ is an invariant direction of the $S^1$-action on
   $\Hom_{{\cal O}_C}
        (\,{\cal I}_{\alpha,\beta}\,,\,{\cal O}_C\,)$,
   then $t\cdot f = t^{\mu} f$ for some $\mu\in{\Bbb Z}$.
   From the above expression, this means that
   $$
    (\,t^{\alpha}\,h_0(t^{-1}z)\,,\,
                   t^{-\beta}\,h_{\infty}(tw)\,)\;
    =\; (\,t^{\mu}\,h_0(z)\,,\, t^{\mu}\,h_{\infty}(w)\,)\;
    \hspace{1ex}\mbox{for all $t$}\,,
   $$
   which implies that
   $$
    f\;=\; f_\mu\;
    :=\; (\,h_0(z)\,,\,h_{\infty}(w)\,)\;
     =\; (\,c\,z^{\alpha-\mu}\,,\, c\,w^{\beta+\mu}\,)\,.
   $$
  From this, one concludes that
   $$
    -\beta\;\le\;\mu\;\le\;\alpha\,,\hspace{1em}
     \mu\,\in\,{\Bbb Z}\,;
   $$
   with the associated weight subspace spanned by $f_{\mu}$.

 \medskip
 \item []
  \hspace{-1em}(2.7)\hspace{1ex}
  Resume the indices $(i,i^{\prime})$
   for ${\cal I}_{\alpha_i,\beta_{i^{\prime}}}$.
   Then

 \medskip
 \item []
  {\bf Lemma [weight subsystem $\Wt_{\,3}$].} {\it
  \begin{itemize}
   \item [{\rm (1)}]
    Let $\Wt_{\,3}^{\;\prime}$ be the system of weights of
     the $S^1$-action on
     $\Hom_{{\cal O}_C}(\,{\cal V}\,,\,{\cal O}_C\,)$.
    Then the weight system $\Wt_{\,3}$ for the $S^1$-action on
    $$
     \Hom_{{\cal O}_C}(\,{\cal V}\,,\,
      {\cal E}^n\!/
        \mbox{\raisebox{-.4ex}{$\widehat{\cal V}$}}\,)\;
    $$
    is given by $\Wt_{\,3}=(n-r)\,\Wt_{\,3}^{\;\prime}$,
    i.e.\ same set of integers ${\mu}$ as in $\Wt_{\,3}^{\;\prime}$ but
    with multilicity $m_{\mu}=(n-r)\,m_{\mu}^{\prime}$, 

  \item [{\rm (2)}]
   $\Wt_{\,3}^{\;\prime}$ is given by 
   $$
    \Wt_{\,3}^{\;\prime}\;=\;\bigsqcup_{\,i=1}^{\,r}\,
      (\,[\,-\beta_{i^{\prime}}\,,\,\alpha_i\,]\,\cap\,{\Bbb Z}\,)\,.
   $$
   Recall
    $(\alpha_1,\,\ldots\,,\alpha_r,\,;\,\beta_1,\,\ldots\,,\beta_r)$
   rewritten as
    $$
     \begin{array}{ccccccccccccccc}
      0 & \le  & a_1 & <  & \cdots & < & a_k\,(=\alpha_r)\,;
        & &   0 & \le  & b_1 & <  & \cdots & < & b_l\,(=\beta_r)
                                                        \\[.6ex]
        &      & m_1 &    & \cdots &   & m_k
        & &     &      & n_1 &    & \cdots &   & n_l
     \end{array}
    $$
    with the multiplicity indicated.
    Then any $\nu\in\,[\,-\beta_r\,,\,\alpha_r\,]\cap{\Bbb Z}$
    is in $\Wt_{\,3}^{\;\prime}$. Its multiplicity $m_{\mu}^{\prime}$
    in $\Wt_{\,3}^{\;\prime}$ is given by
    $$
     m_{\mu}^{\prime}\;
      =\;\left\{
         \begin{array}{ll}
           n_l\,+\, \cdots\, + n_{j}
            & \mbox{if $\; -b_j\;\le\;\mu\;<\;-b_{j-1}\,.$} \\[.6ex]
           r
            & \mbox{if $\; -b_1\;\le\;\mu\;\le\;a_1\,,$}   \\[.6ex]
           m_k\,+\, \cdots\, + m_j
            & \mbox{if $\; a_{j-1}\;<\;\mu\;\le\;a_j\,,$}
         \end{array}
        \right.
    $$
  \end{itemize}
  } 

 \medskip
 \item []
  From this expression, it is clear that $\Wt_{\,3}$ depends only on
   $(\alpha_1,\,\ldots,\alpha_r\,;\,\beta_1,\,\ldots\,\beta_r)$
   and hence only on the connected component of
   the $S^1$-fixed-point locus, as it should.

 \medskip
 \item []
 {\it Proof of Lemma.}
  Consider the two sets of lattice points in
   ${\Bbb Z}\oplus{\Bbb Z}\subset{\Bbb R}^2\,$:
   $$
    A\;=\; \{\,(\alpha_i, r-i+1)\,|\, i=1,\,\ldots,\, r\,\}
     \hspace{1em}\mbox{and}\hspace{1em}
    B\;=\; \{\,(-\beta_i, r-i+1)\,|\, i=1,\,\ldots,\, r\,\}\,,
   $$
   and the $r$-many line segments connecting
   $(-\beta_{i^{\prime}},r-i^{\prime}+1)$ and $(\alpha_i,r-i+1)$.
   Let $\pi$ be the projection of ${\Bbb R}^2$ to the horizontal
   axis $L\supset{\Bbb Z}$.
  Then, for an integer $\mu\in{\Bbb Z}\subset L$,
   the multiplicity $m_{\mu}$ of $\mu$ in $\Wt_{\,3}^{\;\prime}$
   is the same as the number of the line segments above
   whose projection into $L$ contain $\mu$.
   Thus, $m_{\mu}>0$ if and only if $\mu\in[-\beta_r, \alpha_r]$.
  To read off $m_{\mu}$, one combs the collection of line
   segments so that each line segment becomes a three-edged-path
   with the first and the third edge horizontal and the middle
   one vertical and contained in the vertial axis,
  cf.\ {\sc Figure} 2-2-1.
  \begin{figure}[htbp]
   \setcaption{{\sc Figure} 2-2-1.
    \baselineskip 14pt
    The mutiplicity of $\mu$ and the combing of the line segments.
   } 
   \centerline{\psfig{figure=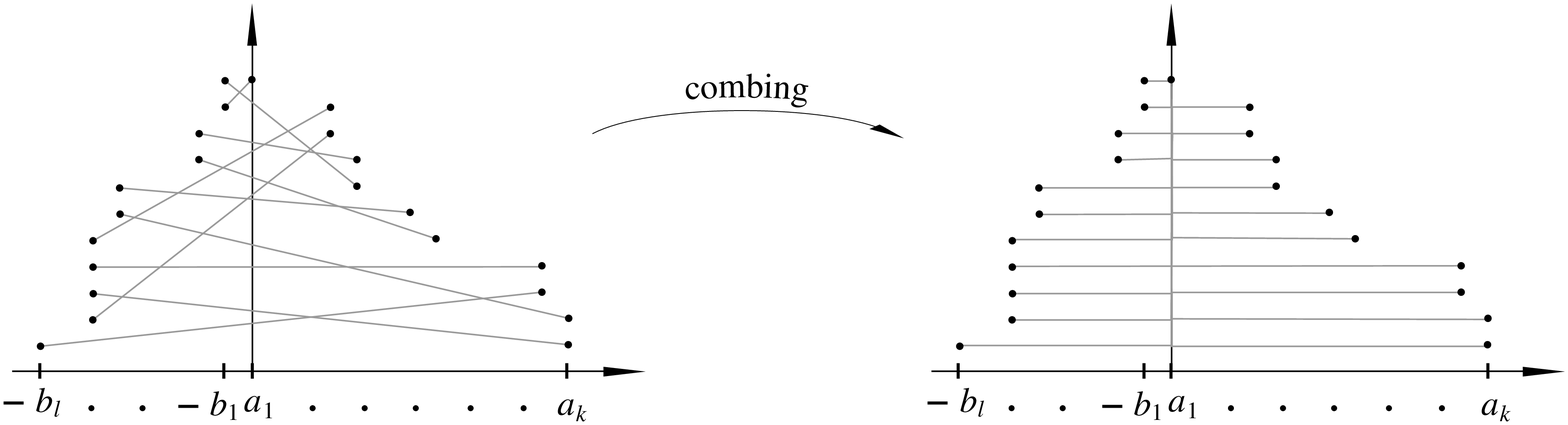,width=13cm,caption=}}
  \end{figure}
  From this, one concludes $m_{\mu}$ as stated in the Lemma.
  This concludes the proof.
\end{itemize}
\noindent\hspace{14cm}$\Box$

\bigskip

To summarize$\,$:

\bigskip

\bigskip 

\noindent 
{\bf Theorem 2.2.1 [$S^1$-weight].} {\it
 The $S^1$-weights on the tangent space
 $$
  T_{({\cal E}^n\rightarrow
   {\cal E}^n\!/\mbox{\raisebox{-.4ex}{\scriptsize${\cal V}$}}
   )}\Quot_P({\cal E}^n)\;
  \simeq\; \Hom_{{\cal O}_C}(\,
            {\cal V}\,,\,
            {\cal E}^n\!/\mbox{\raisebox{-.4ex}{${\cal V}$}}\,)
 $$
 of Quot-scheme $\Quot({\cal E}^n)$ at an $S^1$-fixed-point
 $({\cal E}^n\rightarrow
   {\cal E}^n\!/\mbox{\raisebox{-.4ex}{${\cal V}$}})$
 are the disjoint union of $\Wt_{\,1}$, $\Wt_{\,2}$, and $\Wt_{\,3}$, 
 as given above.
} 

\bigskip

This concludes the computation of the weight system.
We now turn to the combinatorics of this system.

\bigskip

\bigskip

\subsection{Combinatorics of the $S^1$-weight system
            and the multiplicity of $0$.}

A generating function for the multiplicity of weights in $\Wt_{\,3}$
 is immediate, following same argument as in the counting of the states
 at various levels in conformal field theory, e.g.\ [G-S-W].
 An example is given by the following formal function
 $$
  \prod_{j=0}^{\infty}\,
        \frac{1}{1-q_0^{n-r}\cdots q_j^{n-r}\,s^j\,t}\,.
 $$

It remains unclear to us whether the weight systems $\Wt_{\,1}$ and
 $\Wt_{\,2}$ also have elegant generating functions; nevertheless
 they can be obtained from the following manipulations.

\bigskip

\begin{flushleft}
{{\it $\bullet$ The weight subsystem $\Wt_{\,1}$}$\,$:}
\end{flushleft}
\begin{itemize}
 \item [(1)]
  Consider the formal expansion
  $$
    \left.\left(\,
     \prod_{j=0}^{\infty}\,
           \frac{1}{1-q_{-j}\,\cdots\,q_0\,s^j\,t\,v_j}\,
                 \right)\!\right|_{\,q_0=1}\;
      =\; \sum_{k,\,l,\,P}\, A^{(1)}_{k,\,l,\,P}({\bf q})\:
                             A^{(2)}_P({\bf v})\:s^k\,t^{\,l}\,,
  $$
  where ${\bf q}$ and ${\bf v}$ represent collectively
  the two sets of variables $q_i$ and $v_i$ respectively.
  Note that both $A^{(1)}_{k,\,l,\,P}({\bf q})$ and
  $A^{(2)}_P({\bf v})$ are monomials.

 \item [(2)]
  Do the substitutions
  $$
   A^{(2)}_P({\bf v})\;
    \stackrel{
      v_j\rightarrow A^{(1)}_{k,\,l,\,P}({\bf q}_{i\rightarrow i+j+1})
       \rule[-1ex]{0ex}{1ex}
             }{\mbox{------------$\rightarrow$}}\;
     \widehat{A^{(2)}_P}({\bf q})\,,
  $$
  where ${\bf q}_{i\rightarrow i+j+1}$ means that
   $q_i$ is replaced by $q_{i+j+1}$ for all $i$.
  The result $\widehat{A^{(2)}_P}({\bf q})$ is a monomial in
   ${\bf q}$ and the multiplicity of $j\in{\Bbb Z}$ is $n_j$
   if $q_j^{n_j}$ appears as a primary factor of
   $\widehat{A^{(2)}_P}({\bf q})$.
\end{itemize}
(Cf.\ See Example 2.3.2 below.)

\bigskip

\begin{flushleft}
{{\it $\bullet$ The weight subsystem $\Wt_{\,2}$}$\,$:}
\end{flushleft}
\begin{itemize}
 \item [(1)]
  Consider the formal expansion
  $$
    \left.\left(\,
     \prod_{j=0}^{\infty}\,
           \frac{1}{1-q_0\,\cdots\,q_j\,s^j\,t\,v_j}\,
                 \right)\!\right|_{\,q_0=1}\;
     =\; \sum_{k,\,l,\,P}\, B^{(1)}_{k,\,l,\,P}({\bf q})\:
                            B^{(2)}_P({\bf v})\:s^k\,t^{\,l}\,.
  $$

 \item [(2)]
  Do the substitutions
  $$
   B^{(2)}_P({\bf v})\;
    \stackrel{
      v_j\rightarrow B^{(1)}_{k,\,l,\,P}({\bf q}_{i\rightarrow i-j-1})
       \rule[-1ex]{0ex}{1ex}
             }{\mbox{------------$\rightarrow$}}\;
     \widehat{B^{(2)}_P}({\bf q})\,,
  $$
  where ${\bf q}_{i\rightarrow i-j-1}$ means that
   $q_i$ is replaced by $q_{i-j-1}$ for all $i$.
  The result $\widehat{B^{(2)}_P}({\bf q})$ is a monomial in
   ${\bf q}$ and the multiplicity of $j\in{\Bbb Z}$ is $n_j$
   if $q_j^{n_j}$ appears as a primary factor of
   $\widehat{B^{(2)}_P}({\bf q})$.
\end{itemize}

\bigskip

\noindent
{\it Remark 2.3.1.}
\begin{itemize}
 \item [(1)]
  The powers $k$ and $l$ and the monomials
   $A^{(1)}_{k,\,l,\,P}({\bf q})$ and $A^{(2)}_P({\bf v})$
   are related as follows.
   $A^{(2)}_P({\bf v})$ is the monomial that encodes the partition
   $P$ of $k$ into the summmation $l$-many non-negative integers.
   Corresponding to $P$ is a conjugate partition $\widehat{P}$.
  The Young diagram associated to $P$ is conjugate to that associated
   to $\widehat{P}$.
  The monomial $A^{(1)}_{k,\,l,\,P}({\bf q})$ is determined
   by the partition $\widehat{P}$.
  Similarly for $B^{(1)}_{k,\,l,\,P}({\bf q})$ and
   $B^{(2)}_P({\bf v})$.
  (Cf.\ Example 2.3.2 below.)

 \item [(2)]
  Item (1) above implies that for a distinguished $S^1$-fixed-point
   component $F_{\alpha_1,\,\ldots,\,\alpha_r\,;\,0,\,\ldots,\,0}$,
   the subweight system $\Wt_{\,1}$ is generated completely by
   the Young diagram associated to $d=\alpha_1+\,\cdots\,+\alpha_r$
   as a partition of $d$ by a ``partial tensor'' with the conjugate
   Young diagram, as illustrated in {\sc Figure} 2-3-1.
\end{itemize}
\begin{figure}[htbp]
 \setcaption{{\sc Figure} 2-3-1.
  \baselineskip 14pt
   Generation of $\Wt_{\,1}$ from a single Young diagram.
   In the final diagram, the vertical scale is only $1/4$ of
    the horizontal scale.
 } 
 \centerline{\psfig{figure=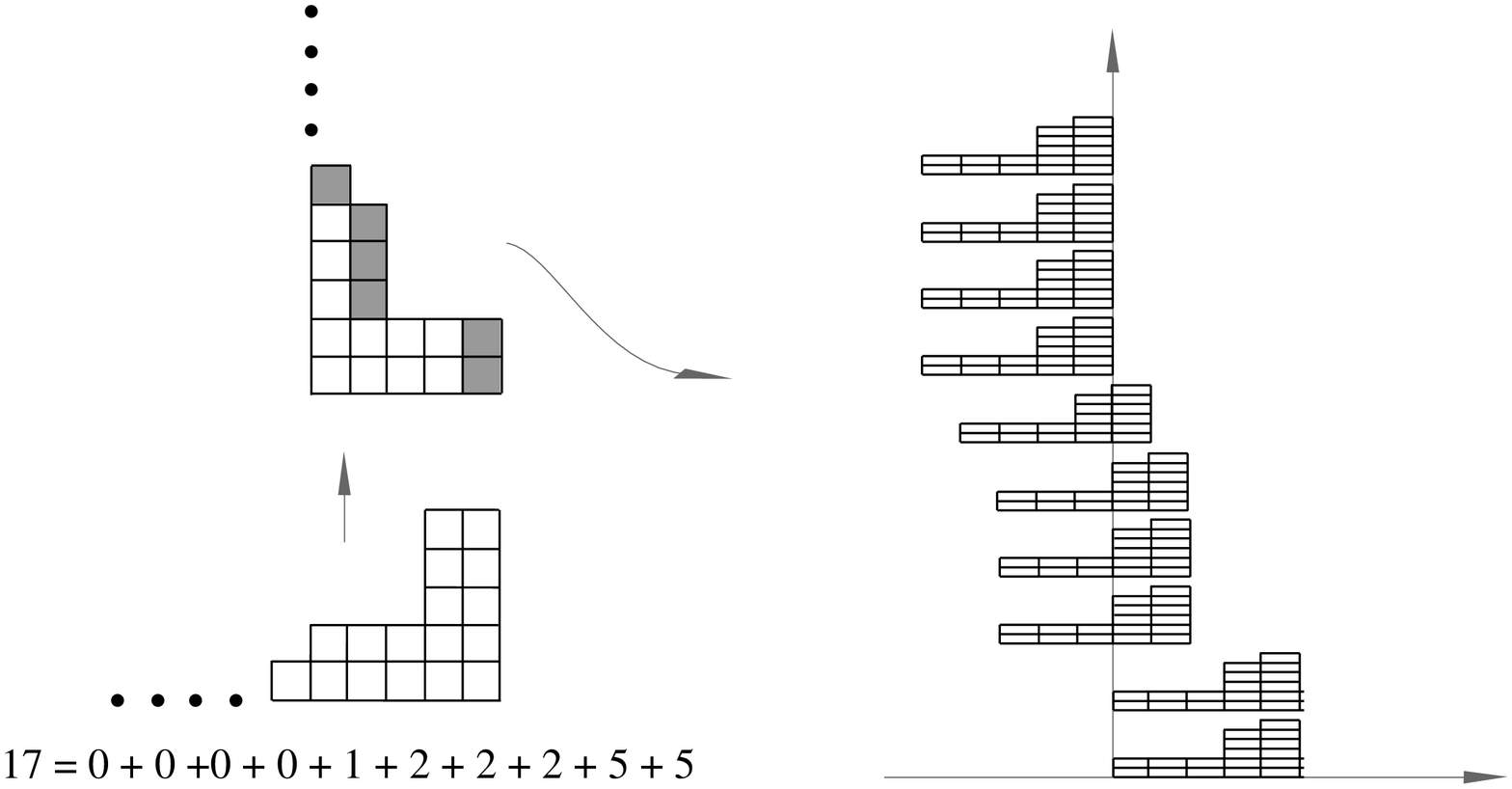,width=13cm,caption=}}
\end{figure}

\bigskip

\noindent
{\bf Example 2.3.2 [weight computation for $\Wt_{\,1}$].}
Consider the expansion
$$
 \left.\left(\,
  \prod_{j=0}^{\infty}\,
        \frac{1}{1-q_{-j}\,\cdots\,q_0\,s^j\,t\,v_j}\,
              \right)\!\right|_{\,q_0=1}\;
  =\; \sum_{k,\,l,\,P}\, A^{(1)}_{k,\,l,\,P}({\bf q})\:
                         A^{(2)}_P({\bf v})\:s^k\,t^{\,l}\,.
$$
Consider, for example, the case $r=10$ and $d=17$.
Then, to determine the weight subsystem $\Wt_{\,1}$ for
 the normal bundle
 to the component $F_{0,0,0,0,1,2,2,2,5,5\,;\,0,\,\cdots,\,0}$ in
 $E_0$, one only needs to look at the (unique) term in the expansion
 with $A^{(2)}_P({\bf v})\,=\,v_0^4\,v_1\,v_2^3\,v_5^2$,
 corresponding to the partition $P\,:\,17=0+0+0+0+1+2+2+2+5+5\,$:
 \begin{eqnarray*}
  \lefteqn{
   (\,A^{(1)}_{17,\,10,\,P}({\bf q})\,A^{(2)}_P({\bf v})\,
                             s^{17}\,t^{10}\,)|_{q_0=1} }  \\[.6ex]
   & & =\; q_{-1}\,(q_{-2}\,q_{-1})^3\,
             (q_{-5}\,q_{-4}\,q_{-3}\,q_{-2}\,q_{-1})^2\:
            v_0^4\,v_1\,v_2^3\,v_5^2\:s^{17}\,t^{10}  \\[.6ex]
   & & =\; q_{-5}^2\,q_{-4}^2\,q_{-3}^2\,q_{-2}^5\,q_{-1}^6\:
             v_0^4\,v_1\,v_2^3\,v_5^2\:s^{17}\,t^{10}\,.
 \end{eqnarray*}
 Observe that the conjugate partition $\widehat{P}\,$:
  $17=0+0+0+0+0+2+2+2+5+6$ is encoded in the monomial in $q$'s.
 Now do the substitiution with the rule of shifting the indices
 as given above$\,$:
 \begin{eqnarray*}
  \lefteqn{
    v_0^4\,v_1^1\,v_2^3\,v_5^2\;\longrightarrow }    \\[.6ex]
    & &
     (\,q_{-4}^2\,q_{-3}^2\,q_{-2}^2\,q_{-1}^5\,q_{\bf \,0}^6\,)^{\bf 4}\:
     (\,q_{-3}^2\,q_{-2}^2\,q_{-1}^2\,q_0^5\,q_{\bf 1}^6\,)^{\bf 1}\:
     (\,q_{-2}^2\,q_{-1}^2\,q_0^2\,q_1^5\,q_{\bf 2}^6\,)^{\bf 3}\:
     (\,q_1^2\,q_2^2\,q_3^2\,q_4^5\,q_{\bf 5}^6\,)^{\bf 2}\,
     |_{\,q_0=1}\; \\[.6ex]
    & &
    =\; q_{-4}^8\,q_{-3}^{10}\,q_{-2}^{16}\,q_{-1}^{28}\,
        q_1^{25}\,q_2^{22}\,q_3^4\,q_4^{10}\,q_5^{12}\,,
 \end{eqnarray*}
 where some of the indices are boldfaced to make the pattern manifest.
 Let $\alpha$ be the generator of $H^{\ast}_{S^1}(\pt)$. Then
 the $S^1$-weights in $\Wt_{\,1}$ for the normal bundle to
 $F_{0,0,0,0,1,2,2,2,5,5\,;\,0,\,\cdots,\,0}$ is
 $$
  8(-4\alpha),\, 10(-3\alpha),\, 16(-2\alpha),\, 28(-\alpha),\,
   25(\alpha),\,  22(2\alpha),\,   4(3\alpha),\, 10(4\alpha),\,
  12(5\alpha)\,.
 $$

\noindent\hspace{14cm}$\Box$

\bigskip

Though it can be obtained also from Sec.\ 2.2, the following lemma
follows immediately from the combinatorics of the weight system
discussed in this subsection.

\bigskip

\noindent
{\bf Lemma 2.3.3 [multiplicity of $0$].} {\it
 The multiplicity of $0$ in the $S^1$-weight system $\Wt$ to
  the restriction of the tangent bundle $T_{\ast}\Quot_P({\cal E}^n)$
  to $F_{\alpha_1,\,\ldots,\alpha_r\,;\,\beta_1,\,\ldots,\,\beta_r}$
  is equal to
  $\dimm F_{\alpha_1,\,\ldots,\alpha_r\,;\,\beta_1,\,\ldots,\beta_r}$.
 Consequently, the $S^1$-weight system of the normal bundle
  to $F_{\alpha_1,\,\ldots,\alpha_r\,;\,\beta_1,\,\ldots,\beta_r}$
  is exactly the subsystem of non-zero weights in $\Wt$.
} 

\bigskip

\noindent
{\bf Corollary 2.3.4 [$\,e_{S^1}$ invertible].} {\it
 Let $E$ be any of the $S^1$-fixed-point component in
  $\Quot_P({\cal E}^n)$.
 Then the $S^1$-weights of the normal bundle
  $\nu_E(\Quot_P{\cal E}^n)$ to $E$ are all nonzero.
 Consequently $e_{S^1}(\nu/\Quot_P({\cal E}^n))$ is invertible
  in $A^{\ast}(E)(\alpha)$, where $\alpha$ is a generator of
  the ring $H^{\ast}_{S^1}(\pt)$.
} 

\bigskip

\noindent
{\it Proof of Lemma 2.3.3.}
 Recall the three subsystems $\Wt=\Wt_1+\Wt_2+\Wt_3$ from Sec.\ 2.2.
 The multiplicity of $0$ in $\Wt_{\,3}$ is $(n-r)r$.
 For the weight subsystem $\Wt_{\,1}$, from the above discussion
 on the Young tableau associated to $\Wt_{\,1}$ and also 
 the characteristic function $\chi_m^A$ for $\Wt_{\,1}$
 defined in Sec.\ 2.2, one has that the multiplicity of $0$
 in $\Wt_{\,1}$ is given by
 $$
  m_{k-1}m_k+m_{k-2}(m_k+m_{k-1})+\,\cdots\,
    +\, m_1(m_l+\cdots+m_2)\,.
 $$
 Similarly the multiplicity of $0$ in the weight subsystem
 $\Wt_{\,2}$ is given by
 $$
  n_{l-1}n_l+n_{l-2}(n_l+n_{l-1})+\,\cdots\,
    +\, n_1(n_l+\cdots+n_2)\,.
 $$
 Consequently the multiplicity of $0$ in $\Wt$ is given by
 \begin{eqnarray*}
  \lefteqn{
    m_{k-1}m_l+m_{k-2}(m_k+m_{k-1})+\,\cdots\,
                           +\, m_1(m_l+\cdots+m_2)  } \\[.6ex]
    & & +\, n_{l-1}n_l+n_{l-2}(n_l+n_{l-1})+\,\cdots\,
                                +\, n_1(n_l+\cdots+n_2)\,
    +\, (n-r)r\,.
 \end{eqnarray*}

 On the other hand,
 {\footnotesize
 \begin{eqnarray*}
  \lefteqn{
   \dimm F_{\alpha_1,\,\ldots,\alpha_r\,;\,
            \beta_1,\,\ldots\,,\beta} } \\[.6ex]
   & &
     =\;\dimm(\Fl_{m_1,\,m_1+m_2,\,\ldots,\,m_1+\cdots m_{k-1},\,r}
                                         ({\footnotesizeBbb C}^n))\,
       +\,
       \dimm(\Fl_{n_1,\,n_1+n_2,\,\ldots,\,n_1+\cdots +n_{l-1},\,r}
                                         ({\footnotesizeBbb C}^n))\,
       -\, \dimm (\Gr_r({\footnotesizeBbb C}^n))   \\[.6ex]
    & &  =\; (n-r)r\,
           +\, (m_1+\cdots+m_{k-1})m_k\, +\,
                   (m_1+\cdots+m_{k-2})m_{k-1}\,+\, \cdots\,
              +\, (m_1+m_2)m_3\, +\, m_1m_2   \\[.6ex]
    & &  \hspace{2em}
           +\, (n-r)r\,
               +\,(n_1+\cdots+n_{l-1})n_l\, +\,
                   (n_1+\cdots+n_{l-2})n_{l-1}\,+\, \cdots\,
                +\, (n_1+n_2)n_3\, +\, n_1n_2 \\[.6ex]
    & &  \hspace{2em}
           -\, (n-r)r\,,
 \end{eqnarray*}
 {\normalsize where}} 
 we have used that fact that
 $m_1+\,\cdots\,+m_k\;=\; n_1+\,\cdots\,+ n_l\;=\; r$.
 By rearrangement of terms, we see that this is the same as
 the multiplicity of $0$ and hence conclude the lemma.

\noindent\hspace{14cm} $\Box$

\bigskip

\bigskip

\section{Mirror principle computation for Grassmannian manifolds.}

\bigskip

\subsection{\bf The distinguished $S^1$-fixed-point components
   and the hyperplane-induced class.}

To make the comparison immediate, here we follow the notations
 in [L-L-Y1$\,$: III, Sec.\ 5.4].
Recall the following approach ibidem to compute $A(t)$
 when there is a commutative diagram$\,$:
$$
 \begin{array}{cccccl}
  F_0  & \stackrel{e^Y}{\longrightarrow}  & Y_0
       & \stackrel{g}{\longleftarrow}    & E_0  & \\[.6ex]
  \hspace{1ex}\downarrow\mbox{\scriptsize $i$}
     & & \hspace{1ex}\downarrow\mbox{\scriptsize $j$}
     & & \hspace{1ex}\downarrow\mbox{\scriptsize $k$} &  \\
  M_d  & \stackrel{\varphi}{\longrightarrow}  & W_d
       & \stackrel{\psi}{\longleftarrow} & {\cal Q}_d  & ,
 \end{array}
$$
where
 ${\cal Q}_d$ is an $S^1$-manifold,
 $\psi:{\cal Q}_d\rightarrow W_d$ is an $S^1$-equivariant resolution
  of singularities of $\varphi(M_d)$, 
 $E_0$ is the set of fixed-points in $\psi^{-1}(Y_0)$ and is
  called the distinguished $S^1$-fixed-point component, and
 $\varphi_{\ast}[M_d]=\psi_{\ast}[{\cal Q}_d]$ in $A_{\ast}^{S^1}(W_d)\,$.

In the current case, $X$ is the Grassmannian manifold
 $\Gr_r({\Bbb C}^n)$, ${\cal Q}_d$ is the Quot-scheme
 $\Quot_{P(t)=(n-r)t+(d+n-r)}({\cal E}^n)$,
 and the linear sigma model $W_d$ for $X$ is the projective space
 ${\Bbb P}(H^0(C,\,{\cal O}C(d))\otimes \Lambda^r{\Bbb C}^n)$
 of $(\!\begin{array}{c}
         \mbox{\scriptsize $n$} \\[-1.2ex]
         \mbox{\scriptsize $r$}
        \end{array}\!)$-tuple of degree-$d$ homogeneous polynomials
 on $C$.
 This is a linear sigma model for ${\Bbb P}(\Lambda^r{\Bbb C}^n)$
 that is turned into a linear sigma model for $X$ via the Pl\"{u}cker
 embedding
 $\Gr_r({\Bbb C}^n)\rightarrow{\Bbb P}(\Lambda^r{\Bbb C}^n)$.

An element in $W_d$ can be written as
 $$
  [\,\sum_j c_{1j}z_0^jz_1^{d-j}\,:\,
     \sum_j c_{2j}z_0^jz_1^{d-j}\,:\,
           \,\cdots\,\,]\,,
 $$
 where $[z_0:z_1]$ is the homogeneous coordinates for $C$ and
 $c_{ij}\in{\Bbb C}$ with
  $1\le i\le(\!\begin{array}{c}
                \mbox{\scriptsize $n$} \\[-1.2ex]
                \mbox{\scriptsize $r$}
               \end{array}\!)$
  and $0\le j\le d$.
The group $S^1$ acts on $W_d$ by
 $$
  [\,\sum_j c_{1j}z_0^jz_1^{d-j}\,:\,\cdots\,\,]\;
   \longmapsto\; 
   [\,\sum_j c_{1j}(tz_0)^jz_1^{d-j}\,:
           \,\cdots\,\,]\,, \hspace{1em} t\in S^1\,.
 $$
There are $(d+1)$-many $S^1$-fixed-point components in $W_d$, each
 of which consists of points of the form
 $[\,c_{1j}z_0^jz_1^{d-j}\,:\,c_{2j}z_0^jz_1^{d-j}\,:\,\cdots\,\,]$
 for $0\le j\le d$ and is isomorphic to
 ${\Bbb P}(\Lambda^r({\Bbb C}^n))$.
From [L-L-Y1$\,$: II, III], the $S_1$-fixed-point component $F_0$
 in $M_d$ consists of degree-$(1,d)$ stable maps $(C,f)$ into
 $C\times\Gr_r({\Bbb C}^n)
             \subset C\times{\Bbb P}(\Lambda^r{\Bbb C}^n)$
 that is obtained by gluing a degree-$(1,0)$ stable map
 $(C_1=\CP^1,f_1,\infty)$ and a degree-$(0,d)$ stable map
 $(C_2, f_2, x)$ with $f_1(\infty)=f_2(x)$ at their marked point.
Regard these as stable maps into the projective space
 ${\Bbb P}(\Lambda^r{\Bbb C}^n)$, then
 [L-L-Y1$\,$: I, Sec.2, Example 10 and III, Sec.\ 3] implies that
 the $S^1$-fixed-point component $Y_j$ in $W_d$ consists
 of point of the form
 $[\,c_{1d}z_0^{d-j}z_1^j\,:\,c_{2d}z_0^{d-j}z_1^j\,:\,\cdots\,\,]$.
In particular, $Y_0$ consists of points of the form
 $[\,c_{1d}z_0^d\,:\,c_{2d}z_0^d\,:\,\cdots\,\,]$.

The map
 $$
  \psi\;:\; {\cal Q}_d=\Quot_{P(t)}({\cal E}^n)\;
    \longrightarrow\;
    W_d
    ={\Bbb P}(H^0(C,\,{\cal O}_C(d))
                                  \otimes \Lambda^r{\Bbb C}^n)
 $$
 is given as follows.
Write $C=\Proj{\Bbb C}[z_0,z_1]$, where ${\Bbb C}[z_0,z_1]$ is
 regarded as a graded ring with grading given by the total degree.
 Then ${\cal E}^n$ is the sheaf associated to the graded
 ${\Bbb C}[z_0,z_1]$-module
 ${\frak M}:={\Bbb C}[z_0,z_1]^{\oplus\,n}$,
 whose grade-$d$ piece ${\frak M}_d$ is given by 
 $$
  {\frak M}_d\;=\;\{\,(f_1, \,\cdots,\, f_n)\,|\, f_i \hspace{1ex}
           \mbox{homogeneous polynomial of $d$ in $z_0,z_1$}\,\}\,.
 $$
A point
 $({\cal E}^n
   \rightarrow {\cal E}^n\!/\mbox{\raisebox{-.4ex}{${\cal V}$}})
   \in {\cal Q}_d$
 is the same as a subsheaf ${\cal V}\hookrightarrow{\cal E}^n$,
 which then corresponds to a graded submodule
 ${\frak N}_{\cal V}$ in ${\frak M}$ of rank $r$. 
Let $e_1,\,\ldots,\,e_r\in{\frak M}$ be a basis for
 ${\frak N}_{\cal V}$. Express each $e_i$ as a column vector
 with entries in ${\Bbb C}[z_0,z_1]$ and consider the matrix
 $A_{\cal V}=[\,e_1,\,\ldots,\,e_r\,]$.
When the quotient sheaf 
 ${\cal E}^n\!/\mbox{\raisebox{-.4ex}{${\cal V}$}}$
 has degree $d$, all the $r\times r$-minors of $A_{\cal V}$,
 if not zero, must be of degree $d$ as well.
The map $\psi$ sends
 $({\cal E}^n
   \rightarrow {\cal E}^n\!/\mbox{\raisebox{-.4ex}{${\cal V}$}})$
 then to the
 $(\!\begin{array}{c}
      \mbox{\scriptsize $n$} \\[-1.2ex]
      \mbox{\scriptsize $r$}
     \end{array}\!)$-tuple of $r\times r$-minors of
 $A_{\cal V}$.
(Cf.\ [Ha], [So], [Str], and [S-S].)

Since $\psi$ is $S^1$-equivariant, it sends an $S^1$-fixed-point
 component in ${\cal Q}_d$ into an $S^1$-fixed-point
 component in $W_d$.
To see which $S^1$-fixed-point component in ${\cal Q}_d$ is sent
 to $Y_0$, one only needs to check where a single point in
 $F_{\alpha_1,\,\ldots,\,\alpha_r\,;\,\beta_1,\,\ldots,\,\beta_r}$
 is mapped to.

\bigskip

\noindent
{\bf Lemma 3.1.1.} {\it
 $\psi(F_{\alpha_1,\,\ldots,\,\alpha_r\,;\,
                      \beta_1,\,\ldots,\,\beta_r})\,
  \subset\, Y_{\beta_1+\,\cdots\,+\beta_r}$.
} 

\bigskip

\noindent
{\it Proof.}
Recall that, for a fixed-point
 $({\cal E}^n
   \rightarrow {\cal E}^n\!/\mbox{\raisebox{-.4ex}{${\cal V}$}})
   \in F_{\alpha_1,\,\ldots,\,\alpha_r\,;\,
                         \beta_1,\,\ldots,\,\beta_r}$,
 $$
  \degree{\cal E}^n\!/\mbox{\raisebox{-.4ex}{${\cal V}$}}\;
  =\;\alpha_1+\,\cdots\,+\alpha_r+\beta_1+\,\cdots\,+\beta_r\;
  =\;d\,.
 $$
Observe also that the special fixed-points in
 $F_{\alpha_1,\,\ldots,\,\beta_r}$, for which the two local
 diagonalization match with $\alpha_i\rightarrow\beta_{i^{\prime}}$, 
 corresponds to a subsheaf ${\cal V}$ in ${\cal E}^n$ is isomorphic
 to the direct sum
 $\oplus_i\,{\cal I}_{\alpha_i(0)+\beta_{i^{\prime}}(\infty) }$
 of ideal sheaves ${\cal I}_{\alpha_i(0)+\beta_{i^{\prime}}(\infty)}$
 in ${\cal O}_C$ associated to the degree-$d$ divisor/subscheme
 $\alpha_i(0)+\beta_{i^{\prime}}(\infty)$ in $C$.
Its associated matrix $A_{\cal V}$ can be written as
 $$
  \left(\begin{array}{ccc}
          z_0^{\alpha_1}z_1^{\beta_{1^{\prime}}} & & \\
          & \ddots  &  \\
          & & z_0^{\alpha_r}z_1^{\beta_{r^{\prime}}}\\[.6ex]
          0 & \cdots & \hspace{-1.6ex}0  \\
          \vdots  & \cdots & \hspace{-1.6ex}\vdots \\
          0 & \cdots  & \hspace{-1.6ex}0
        \end{array} \right)
      \hspace{1em}\mbox{with zero entries $a_{ij}$ for $i\ne j$,}
 $$
 after a constant re-trivialization of ${\cal E}^n$.
The $r\times r$-minors of this matrix are all zero except
 the one from the top $r\times r$-submatrix, whose value is
 $z_0^{\alpha_1+\,\cdots\,+\alpha_r}\,
                        z_1^{\beta_1+\,\cdots\,+\beta_r}$.
Thus, $\psi$ maps such point to some 
$$
 [\,0:\,\cdots\,:0:\,
    z_0^{\alpha_1+\,\cdots\,+\alpha_r}\,
                        z_1^{\beta_1+\,\cdots\,+\beta_r}\,:
    0:\,\cdots\,:0\,]\,,
$$
which lies in $Y_{\beta_1+\,\cdots\,+\beta_r}$.
This proves the lemma.

\noindent\hspace{14cm}$\Box$

\bigskip

Since $0\le\beta_1\le \cdots \le \beta_r$, one concludes that 

\bigskip

\noindent
{\bf Corollary 3.1.2 [distinguished components].} {\it
 The distinguished $S^1$-fixed-point locus $E_0$ is given by 
  $$
   E_0\;=\;\coprod_{ 
             \begin{array}{c} 
              \mbox{\scriptsize
                $0\le\alpha_1\le\cdots\le\alpha_r$ } \\[-.8ex]
              \mbox{\scriptsize $\alpha_1+\cdots+\alpha_r=d$}
             \end{array}
                    }\,
          F_{\alpha_1,\,\cdots,\,\alpha_r\,;\,0,\,\ldots\,,\,0}\;,
  $$
  a disjoint union of flag manifolds determined by the multiplicities
  of entries in $(\alpha_1,\,\ldots,\,\alpha_r)$ with
  $\alpha_1+\cdots+\alpha_r=d$.
 {\rm (}Cf.\ Theorem 2.1.9 {\rm )}
} 

\bigskip

On each distinguished $S^1$-fixed-point component
 $F_{\alpha_1,\,\cdots,\,\alpha_r\,;\,0,\,\ldots\,,\,0}$,
 there is the pulled-back hyperplane class
 $k^{\ast}\psi^{\ast}\kappa=g^{\ast}j^{\ast}\kappa$,
 where $\kappa$ is the hyperplane class on $W_d$.
To see what it is, recall first the multiplicity numbers
 $m_1,\,\ldots,\,m_k$ for $0\le\alpha_1\le\,\cdots\,\le\alpha_r\,$
 and the following fact/definition$\,$:

\bigskip

\noindent
{\bf Fact/Definition 3.1.3 [special Schubert cycle].}
 (Cf.\ [Fu1], also [Gr2] and [Jo].) {\rm
 Recall that, over the flag manifold
  $\Fl=\Fl_{m_1,\,m_1+m_2,\,\ldots,\,r}({\Bbb C}^n)$,
  there is a universal flag of bundles
  $S_1\hookrightarrow S_2 \hookrightarrow \,\cdots\,
                    \hookrightarrow S_{k+1}=\Fl\times{\Bbb C}^n$
  with $\rank S_i=m_1+\,\cdots\,+m_i$.
 Then {\it
  the intersection Chow ring $A^{\ast}(\Fl)$ is generated by
  the Chern classes of the quotient bundles
  $S_i/\mbox{\raisebox{-.4ex}{$S_{i-1}$}}$, $1\le i\le k+1$ and
  $S_0=0$, with relations determined by
  $\prod_{i=1}^{k+1} c\,(S_i/\mbox{\raisebox{-.4ex}{$S_{i-1}$}})\,=\,1$}.
 The Schubert cycles that represent these special generators are
  called {\it special Schubert cycles}.
} 

\bigskip

\noindent
Since
 $F_{\alpha_1,\,\cdots,\,\alpha_r\,;\,0,\,\ldots\,,\,0}
  \simeq \Fl_{m_1, m_1+m_2,\,\ldots,\,r}({\Bbb C}^n)$,
 this gives 
 $A^{\ast}(
  F_{\alpha_1,\,\cdots,\,\alpha_r\,;\,0,\,\ldots\,,\,0})$.
Recall also from Sec.\ 2.1 that points in
 $F_{\alpha_1,\,\cdots,\,\alpha_r\,;\,0,\,\ldots\,,\,0}$
 can be represented by $n\times r$-matrices $B(z)$ with
 coefficients in ${\Bbb C}[z]$.
Since the map $g$ is the Pl\"{u}cher embedding and it sends
 $B(z)$ to the tuple of $r\times r$-minors of $B(1)$
 multiplied by the factor $z^d$, the image
 $g(F_{\alpha_1,\,\cdots,\,\alpha_r\,;\,0,\,\ldots\,,\,0})$
 coincides with the image of the Grassmannian manifold
 $\Gr_r({\Bbb C}^n)$ in $Y_0$ via Pl\"{u}cker embedding
and $g$ is indeed the fibration to the base Grassmannian manifold
 given in Theorem 2.1.9.
 %

Let $S\hookrightarrow \Gr_r({\Bbb C}^n)\times{\Bbb C}^n$
 be the universal rank-$r$ bundle over $\Gr_r({\Bbb C}^n)$.
Then the Pl\"{u}cker embedding in the direct bundle language is
 the section from the projectivization of the tautological bundle map
 $\bigwedge^rS=\det S\,
  \hookrightarrow\,Gr_r({\Bbb C}^n)\times\bigwedge^r{\Bbb C}^n$
 over $\Gr_r({\Bbb C}^n)$ and, hence, the hyperplane class on
 $\CP^{\,{n \choose r}-1}$ is pulled back to the Chern class
 $-c_1(S)$ on $\Gr_r({\Bbb C}^n)$ via the Pl\"{u}cker embedding.
On the other hand, the embedding of $Y_0\simeq\CP^{\,{n \choose r}-1}$
 in $W_d\simeq\CP^{\,{n\choose r}\,d+{n\choose r}-1}$ has degree
 $1$ from [L-L-Y1$\,$: I and II].
Together one concludes that$\,$: 

\bigskip

\noindent
{\bf Corollary 3.1.4 [pulled-back hyperplane class].} {\it
 Let $\kappa$ be the hyperplance class in $W_d$.
  Then, with the notation in Fact/Definition 3.1.3, one has 
  $$
   k^{\ast}\psi^{\ast}\kappa\; =\; g^{\ast}j^{\ast}\kappa\;
    =\; -\,c_1(S_k)
  $$
  on the distinguished $S^1$-fixed-point component
  $F_{\alpha_1,\,\cdots,\,\alpha_r\,;\,0,\,\ldots\,,\,0}$.
 Since $S^1$ and ${\Bbb C}^{\times}$ act on these components
  trivially, these classes lift naturally as to classes on
  $(F_{\alpha_1,\,\cdots,\,\alpha_r\,;\,0,\,\ldots\,,\,0})
                                     _{{\tinyBbb C}^{\times}}$
  and will be denoted by the same notation.
} 

\bigskip

\noindent
{\it Remark 3.1.5
 {\rm [}pulled-back hyperplane in Chern roots$\,${\rm ]}.}
 In terms of Chern roots to be discussed in Sec.\ 3.3,
 this class is represented by
 $-\,(y_1+\,\cdots\,+y_r)=y_{r+1}+\,\cdots\,+y_n$.

\bigskip

\bigskip

\subsection{The weight subspace decomposition of the normal bundle
            to the distinguished components.}

In this subsection, we work out an ingredient needed for
the computation of the ${\Bbb C}^{\times}$-equivariant Euler class
of the normal bundle to a distinguished $S^1$-fixed-point component
in Quot-scheme.

\bigskip

\begin{flushleft}
{\bf Reduction of structure group and
     the $S^1$-weight subspaces in matrix forms.}
\end{flushleft}
Note that the notation $P$ in this section is for parabolic subgroups.
Recall that the $\GL(n,{\Bbb C})$-action on ${\Bbb C}^n$ induces
 a $\GL(n,{\Bbb C})$-action on the set of local sections
 in ${\cal E}^n$.
Thus, given a $g\in\GL(n,{\Bbb C})$, one has a correspondence
 ${\cal V}\longmapsto g\cdot{\cal V}$ with a specified isomorphism
 from ${\cal V}$ to $g\cdot{\cal V}$.
This induces a $\GL(n,{\Bbb C})$-action on $\Quot_{P(t)}({\cal E}^n)$,
 which leaves all the $S^1$-fixed-point component invariant.
This $\GL(n,{\Bbb C})$-action on $\Quot_{P(t)}({\cal E}^n)$ commutes
 with the $S^1$-action discussed earlier.
In this way, the normal bundle $\nu$ to a $S^1$-fixed-point
 component $E$ is realized as a homogeneous $\GL(n,{\Bbb C})$-bundle
 and its structure group is the stabilizer $P$ of a point $p$ in
 that component$\,$:
 $\nu_E\,{\cal Q}_d=\GL(n,{\Bbb C})\times_P{\Bbb C}^R$,
 where $R$ is the codimension of $E$ in ${\cal Q}_d$, ${\Bbb C}^R$
 is identified with the fiber of $\nu_E\,{\cal Q}_d$ at $p$ with
 the $P$-action induced from $\GL(n,{\Bbb C})$.

The existence of a flag manifold also as a compact quotient
 implies that one can choose a compact $U(n)$ in $\GL(n,{\Bbb C})$
 such that each $S^1$-fixed-point component is also a $U(n)$-orbit. 
Then the new stabilizer at a point becomes
 $$
  P_0\,=\,U(n)\cap P\;
   =\;U(m_1)\times\,\cdots\,\times\,U(m_k)\times U(n-r)
 $$
 and 
 $$
  \nu_E\,{\cal Q}_d\;
   =\; \GL(n,{\Bbb C})\times_P{\Bbb C}^R\;
   =\; U(n)\times_{P_0}{\Bbb C}^R\,.
 $$
In this way, we have reduced the structure group of
 $\nu_E\,{\cal Q}_d$ to $P_0$ that remains compatible with
 the $S^1$-action.
Applying this to each of the distinguished $S^1$-fixed-point
 components $F_{\alpha_1,\,\ldots,\,\alpha_r\,;\,0,\,\ldots,\,0}$,
 we then realize
 $T_{\ast}\Quot_{P(t)}({\cal E}^n)
     |_{F_{\alpha_1,\,\ldots,\,\alpha_r\,;\,0,\,\ldots,\,0}}$
 as a homogeneous $U(n,{\Bbb C})$-bundle, determined by
 a representation of $P_0$.

Given $0\le\alpha_1\le\,\ldots\,\le\alpha_r$ rewritten as
$$
 \begin{array}{ccccccccccccccc}
  0 & \le   & a_1  & <   & \cdots & < & a_k\,(=\alpha_r)
                                                     \\[.6ex]
    &       & m_1  &     & \cdots &   & m_k
 \end{array}
$$
with the multiplicity indicated, fix a point on
$F_{\alpha_1,\,\ldots,\,\alpha_r\,;\,0,\,\ldots,\,0}$ represented
by the subsheaf ${\cal V}$ in ${\cal E}^n$ determined by
$$
 {\cal V}(U_0)={\Bbb C}[z]\cdot z^{\alpha_1}\oplus\cdots
   \oplus {\Bbb C}[z]\cdot z^{\alpha_r}\oplus 0^{\oplus (n-r)}
 \hspace{1em}\mbox{and}\hspace{1em}
 {\cal V}(U_{\infty})={\Bbb C}[w]^{\oplus r}\oplus 0^{\oplus (n-r)}
$$
(or equivalently, the graded submodule in ${\frak M}$ generated by
   $(0,\,\ldots,\,0,\,z_0^{\alpha_i},\,0\,\ldots,\,0)$ for
   $1\le i\le r$, in the notation of Sec.\ 2.1).
Then $P$ is the subgroup of appropriate block upper triangular
 matrices in $\GL(n,{\Bbb C})$.
Fix a Hermitian inner product on ${\Bbb C}^n$, which renders
 ${\cal E}^n$ a trivialized Hermitain vector bundle, and let
 $U(n)\hookrightarrow\GL(n,{\Bbb C})$ be the subgroup of
 $\GL(n,{\Bbb C})$ with respect to this inner product.
Then the induced action of $U(n)$ on
 $F_{\alpha_1,\,\ldots,\,\alpha_r\,;\,0,\,\ldots,0}$
 is transitive with
 $P_0=P\cap U(n)=U(m_1)\times\,\cdots\,\times U(m_k)\times U(n-r)$
 being the subgroup of $U(n)$ that consists of
 $m_1\times m_1,\,\ldots,\,m_k\times m_k,\, (n-r)\times(n-r)$
 unitary diagonal blocks.

There is an embedding of $\Hom$-groups
 $$
  \Hom_{{\cal O}_C}({\cal V},
            {\cal E}^n\!/\mbox{\raisebox{-.4ex}{${\cal V}$}})
    \hookrightarrow
  \Hom_{{\scriptsizeBbb C}[z]}(
    {\Bbb C}[z]\cdot z^{\alpha_1}\oplus\,\cdots\,
                 \oplus {\Bbb C}[z]\cdot z^{\alpha_r}\,,\,
    {\Bbb C}[z]\cdot \overline{e}_1\oplus\,\cdots\,
                 \oplus {\Bbb C}[z]\cdot\overline{e}_r
                                \oplus {\Bbb C}[z]^{\oplus (n-r)})\,,
 $$
 where the annihilator $\Ann(\overline{e}_i)$ of $\overline{e}_i$
 is the ideal $(z^{\alpha_i})$ in ${\Bbb C}[z]$.
Let $m_0$ be the multiplicity of $0$ in
 $\alpha_1,\,\ldots,\,\alpha_r$.
 Then $\overline{e}_1=\,\cdots\,= \overline{e}_{m_0}=0$
 and the image is a submodule of the latter that consists of
 matrices of polynomials with degree bounds$\,$:
 $$
  f\,=\, [\,f_{ij}(z)\,]_{(n-m_0)\times r}\,,
 $$
 where
 $$
  \degree f_{ij}(z)\;
    \le\; \left\{
      \begin{array}{lcl}
        \alpha_{m_0+i}-1    & \mbox{for}
         & 1\le i \le r-m_0 \hspace{1em}\mbox{and}\hspace{1em}
                                          1\le j \le r\,, \\[.6ex]
        \alpha_j      & \mbox{for}
         & r-m_0+1 \le i\le n-m_0
          \hspace{1em}\mbox{and}\hspace{1em} 1 \le j\le r\,.
      \end{array}
          \right.
 $$
 The $P_0$-action on
 $\Hom_{{\cal O}_C}({\cal V},
             {\cal E}^n\!/\mbox{\raisebox{-.4ex}{${\cal V}$}})$
 is given by
 $$
  f\; \longmapsto\; g\diamond f\,,
    \hspace{1em}\mbox{for
      $f\in \Hom_{{\cal O}_C}({\cal V},
               {\cal E}^n\!/\mbox{\raisebox{-.4ex}{${\cal V}$}})$
       and $g\in P_0$}
 $$
 with 
 $$
  g\diamond f\;
  :=\; g\,\odot\,f\,
        \odot\,\Diag\{z^{-\alpha_1},\,\cdots,\,z^{-\alpha_r}\}\,
        \cdot\,g^{-1}\,\cdot\,
        \Diag\{z^{\alpha_1},\,\cdots,\,z^{\alpha_r}\}\,,
 $$
 where
  $g$ in the formula is the lower-right $(n-m_0)\times(n-m_0)$
   submatrix of the defining matrix of $g$ when acting on ${\Bbb C}^n$,
  $g^{-1}$ is the $r\times r$ upper-left submatrix of
   the inverse of the defining matrix for $g$,
  the operation $\;\cdot\;$ is the usual matrix mutiplication, and
  the operation $\;\odot\;$ is the usual matrix multiplication
   followed by truncations of terms in an entry that exceeds
   the degree bound above.
This shows explicitly that the $P_0$-action and the $S^1$-action on
 $\Hom_{{\cal O}_C}({\cal V},
             {\cal E}^n\!/\mbox{\raisebox{-.4ex}{${\cal V}$}})$
 commute.

From the previous discussions on the $S^1$-weight system,
 each monomial in an entry (a Laurent polynomial in $z$) of
 $$
  \widetilde{f}\;
  :=\; f\,\odot\,\Diag\{z^{-\alpha_1},\,\cdots,\,z^{-\alpha_r}\}
 $$
 gives an $S^1$-invariant subspace in
 $\Hom_{{\cal O}_C}({\cal V},
             {\cal E}^n\!/\mbox{\raisebox{-.4ex}{${\cal V}$}})$.
The degree bound for an entry in $\widetilde{f}$ is given by
 $$
  \left\{
   \begin{array}{lcl}
    -\alpha_j\;\le\;\degree \widetilde{f}_{ij}(z)\;
      \le\; \alpha_{m_0+i}-\alpha_j-1
      & \mbox{for}
      & 1\le i \le r-m_0 \hspace{1em}\mbox{and}\hspace{1em}
                                       1\le j \le r\,, \\[.6ex]
    -\alpha_j\;\le\;\degree \widetilde{f}_{ij}(z)\; \le\; 0
      & \mbox{for}
      & r-m_0+1 \le i\le n-m_0
          \hspace{1em}\mbox{and}\hspace{1em} 1 \le j\le r\,.
   \end{array}
  \right.
 $$
Thus one has a decomposition of the $P_0$-module by the $S_1$-weight
 subspaces, each of which is itself a $P_0$-module$\,$:
 $$
  \widetilde{f}\;
   =\; z^{-\alpha_r}\,\widetilde{f}_{(\alpha_r)}\,
        +\, \cdots\, +
       z^{-1}\,\widetilde{f}_{(1)}\,+\,\widetilde{f}_{(0)}\,
        +\, z\,\widetilde{f}_{(-1)}\,
        +\, \cdots\,
        + z^{\alpha_r-\alpha_1-1}\,
                  \widetilde{f}_{(\,-(\alpha_r-\alpha_1-1)\,)}\,,
 $$
 where the $S^1$-weight for $z^{\mu}$-component here is $-\mu$,
  (cf.\ the expression $z^{\alpha_j-\mu^0_{ij}}$
  in the discusion of the $S^1$-weight system $\Wt_{\,1}$).
 (Note that here we are assuming the generic situation,
  in which $\alpha_1<\alpha_r$ and hence $\alpha_r-\alpha_1-1\ge 0$.
  If $\alpha_1=\,\cdots\,=\alpha_r$, then
   $\alpha_r-\alpha_1-1=-1$ and  
   $\widetilde{f}
    = z^{-\alpha_r}\,\widetilde{f}_{(\alpha_r)}\,
        +\, \cdots\, +
      z^{-1}\,\widetilde{f}_{(1)}\,+\,\widetilde{f}_{(0)}\,$.)

\bigskip

\begin{flushleft}
{\bf The $P_0$-module decomposition of $S^1$-weight spaces and
     the $P_0$-weight system.}
\end{flushleft}

\noindent $\bullet$
{\it $S^1$-weight-subspace decomposition.}
\begin{itemize}
 \item [(1)]
  Recall the multiplicity $m_i$, $1\le i\le k$, of the sequence
   $0\le\alpha_1\le\,\cdots\,\le\alpha_r$ and $m_0$
   the multiplicity of $0$ in the sequence.
  Then the matrices $g$, $g^{-1}$, $f$, $\widetilde{f}$ can be
   put into a block form.
  For example, the $(I,J)$-block for $f$ is
    an $m_I\times m_J$ submatrix if $m_0=0$, or
    an $m_{I+1}\times m_J$ submatrix if $m_0>0$, or
    an $r\times m_J$ submatrix
       if $m_0=0$ and $I=k+1$ or if $m_0>0$ and $I=k$.
   (Cf. {\sc Figure} 3-2-1.)

 \item [(2)]
  In terms of the block form, the decomposition of $\widetilde{f}$
   into a summation of matrices with only one non-zero block
   gives the decomposition of
   $\Hom_{{\cal O}_C}({\cal V},
             {\cal E}^n\!/\mbox{\raisebox{-.4ex}{${\cal V}$}})$
   into representations of $P_0$.
  Consequently, each such summand is the representation of the form 
   $\rho_{m_I}\otimes (\rho_{m_J}^{-1})^t
                  =\rho_{m_I}\otimes\overline{\rho}_{m_J}$,
   where
    $\rho_{m_I}$ is the defining representation of $U(m_I)$,
    $(\rho^{-1})^t$ its inverse transpose, which is the same as
     its complex conjugate $\overline{\rho}$.

 \item [(3)]
  This decomposition is compatible with the $S^1$-weight subspace
   decomposition.
  In fact, {\it the block form of the $S^1$-weight summand
   $\widetilde{f}_{(s)}$,
   $-\alpha_r\le s\le \max\{\alpha_r-\alpha_1-1,0\}$,
   is determined by the Young diagram corresponding to the partition
   $d=\alpha_1+\,\cdots\,+\alpha_r$}.
  They are all {\it ``sparse-lower-triangular" block matrices},
   (cf.\ {\sc Figure 3-2-1}).
  These block forms are invariant under the conjugation followed by
   truncations of terms of degree higher than the upper degree bounds
   $$
    \widetilde{f}_{(s)}\; \longmapsto\;
     g\,\odot\,\widetilde{f}_{(s)}\,\odot\,g^{-1}
   $$
   and hence this gives a decomposition of the homogeneous bundle
   into the direct sum of $S^1$-weight homogeneous subbundles.
  In particular, the lower sub-triangular block form of
    $\widetilde{f}_{(0)}$ corresponds to the tangent bundle
    $T_{\ast}F_{\alpha_1,\,\ldots,\,\alpha_r\,;\,0,\,\ldots,\,0}$.
  The dimension is consistent with the computation in Lemma 2.3.3.  
   
 \item [(4)]
  These sparse-lower-triangular block matrices are determined
  by the Young diagram corresponding to the partition
  $d=\alpha_1\,+\,\cdots\,+\alpha_r$.
  The rule from a Young diagram to the sparse-lower-triangular
   block forms can be summarized in three steps$\,$:
  \begin{itemize}
   \item [(4.1)]
    Take the dual of the Young diagram and put 
     the zero-matrix of the same dimension as $\widetilde{f}$
     into the same block form.
     Copy these zero-matrices by multiplying by a weight factor
      $z^{\nu}$ with
      $-\alpha_r\le\nu\le \max\{\,\alpha_r-\alpha_1-1\,,\,0\}\,$.
       
    \item [(4.2)]
     Recall $a_I$ and $m_I$ at the beginning of this subsection.
     Multiply the dual Young diagram horizontally by
      the multiplicity $m_I$ and fill the block forms
      with these multi-strip as indicated in {\sc Figure} 3-2-1,
      beginning with the block form with $z^{-a_I}$-factor.
     This corresponds to the $S^1$-weight system $\Wt_{\,1}$.

    \item [(4.3)]
      Add all these matrices and fill in all the blocks
       in matrices with negative $\nu$ in $z^{\nu}$ such that
       some block above it is already filled.
      For the block form with factor $z^0$, fill in all
       the blocks in the rows lower than the last filled row.
      These additional filling corresponds to the $S^1$-weight
       system $\Wt_{\,3}$.
  \end{itemize}
\end{itemize}   

\begin{figure}[htbp]
 \setcaption{\small {\sc Figure 3-2-1.}
  \baselineskip 13pt
  The simultaneous decomposition of
   $\Hom_{{\cal O}_C}({\cal V},
             {\cal E}^n\!/\mbox{\raisebox{-.4ex}{${\cal V}$}})$
   by weight subspaces of $S^1$ and representations of $P_0$.
   Original entries in the matrix are divided by light lines while
    blocks are divided by dark lines.
   The think dark line divides the upper $(r-m_0)$ rows and the
    lower $(n-r)$ rows.
   All the unshaded blocks are zero.
  Observe how the block forms are all determined by the Young
   diagram - the conjugate Young diagram is horizontally fattened
   by the various multiplicities and then distribute into
   the block forms (cf.\ the blocks with the same dark shades) -.
  } 
  \centerline{\psfig{figure=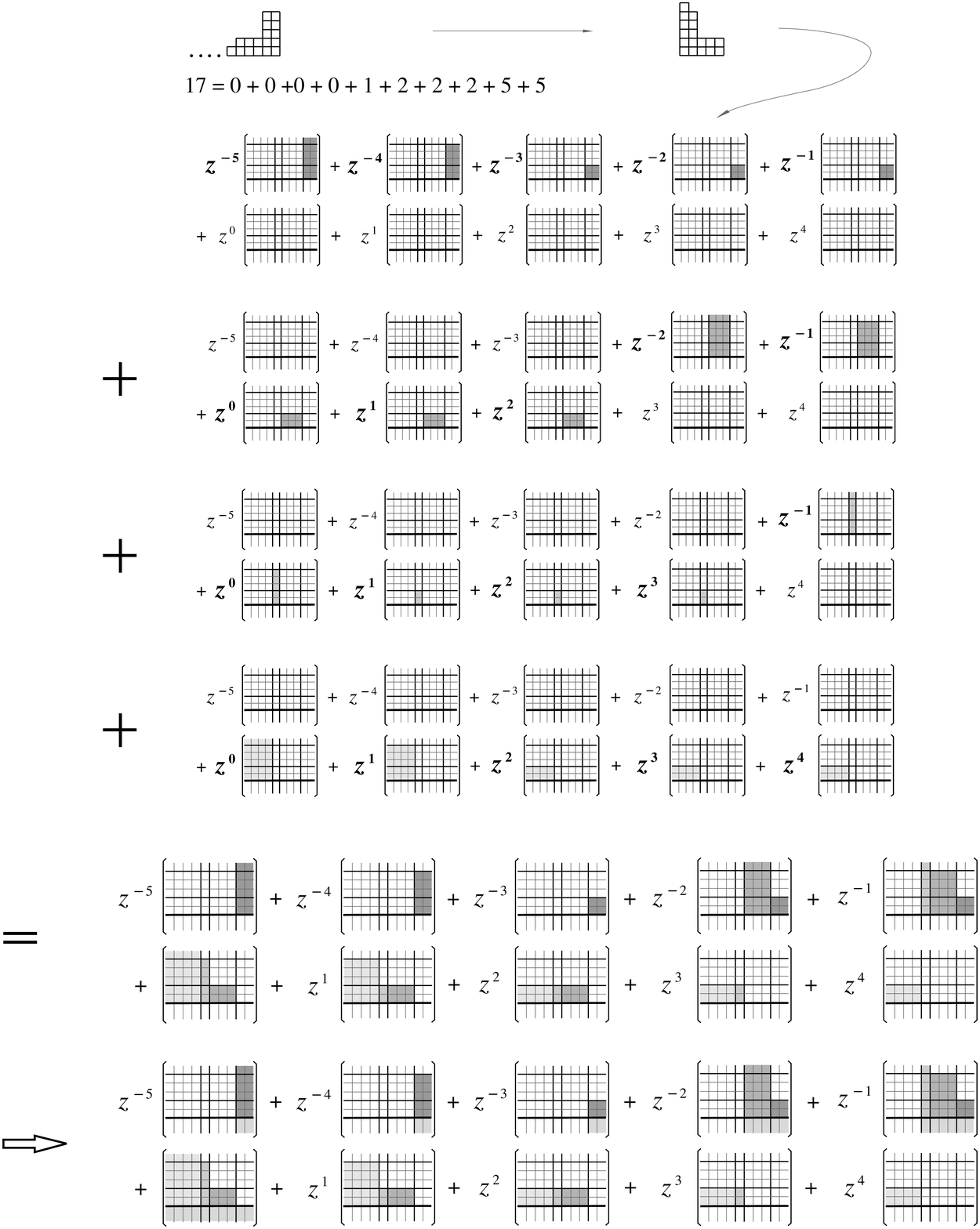,width=14cm,caption=}}
\end{figure}

\begin{itemize}
 \item [(5)]
  The block decomposition of each $S^1$-weight subspace
  into the direct sum of $P_0$-modules.
  \begin{itemize}
   \item [(5.1)]
    {\it Definition of diagonal and off-diagonal blocks$\,$}:
    When $m_0=0$, the diagonal blocks follow the usual definition.
    When $m_0>0$, the diagonal blocks here are the blocks
     that are above and adjacent to the usual diagonal blocks
     (i.e.\ the $(I,I+1)$-blocks).
    All other blocks are called off-diagonal.

   \item [-]
    The diagonal blocks corresponds to a representation
     $\rho_{m_I}\otimes\overline\rho_{m_I}$, where $\rho_{m_I}$
     is the defining representation of some $U(m_I)$.

   \item [-]
    The off-diagonal blocks are irreducible representations
    $\rho_1\otimes\rho_2^{\ast}$ of the product
    $U(m_{I_1})\times U(m_{I_2})$, where $\rho_j$ is the defining
    representation of $U(m_{I_j})$, $i=1,\,2$.

   \item [(5.2)]
    Let
     $(\lambda_1,\,\ldots,\,\lambda_{I_1})$ be the weight system
      of the representation $\rho_1$ and
     $(\lambda^{\prime}_1,\,\ldots,\,\lambda^{\prime}_{I_2})$ 
      be the weight system of the representation $\rho_2$
      (with multiple weight repeated correspondingly),
    then the weight system of $\rho_1\otimes\overline{\rho_2}$
     is given by
     $$
     (\,\lambda_i-\lambda^{\prime}_j\,|\,
                         1\le i\le I_1\,,\; 1\le j\le I_2\,)\,.
     $$
  \end{itemize}
\end{itemize}

\bigskip

\noindent $\bullet$
{\it The $P_0$-weight system $\Wt_{P_0}$ of 
      $\Hom_{{\cal O}_C}({\cal V},
             {\cal E}^n\!/\mbox{\raisebox{-.4ex}{${\cal V}$}})$.}
\begin{itemize}
 \item [(1)]
  Recall the fixed maximal torus the diagonal subgroup
   $T=({\Bbb C}^{\times})^n$ in $P_0$.
  Let $E_{ij}$ be a $(n-m_0)\times r$ matrix with $1$ in
   $(i,j)$-entry and zero elsewhere.
  Then, every subspace of an $S^1$-weight subspace in
   the lower-triangular block that consists of constant multiples
   of some $E_{ij}$ is a $P_0$-weight subspace of weight
   $\lambda_{\,i+m_0}-\lambda_j$.
  Consequently, the $P_0$-weight system can be directly read off
   from the collection of sparse lower-triangular block forms
   obtained from the Young diagram corresponding to the distinguished
   $S^1$-fixed-point component. In expression,
  $$
   \Wt_{P_0}(\,\mbox{Young diagram}\,)\;
     =\;\bigsqcup_{(\,\mbox{\footnotesize
                             triangular block form}\:\Delta\,)}\;
        \bigsqcup_{(\,\mbox{\footnotesize block}\:\Box\in\Delta\,)}\;
        \bigsqcup_{(i,j)\in\Box}\,(\lambda_{\,i+m_0}-\lambda_j)
  $$
  where the Young diagram is the one corresponding to
    $F_{\alpha_1,\,\ldots,\,\alpha_r\,;\,0,\,\ldots,\,0}$
    (namely the partition $d=\alpha_1+\,\cdots\,+\alpha_r$),
   $\sqcup$ means the disjoint combination with multiplicity
     allowed,
   $\Box\in\Delta$ means that the block is in the triangular
     block form $\Delta$, and
   $(i,j)\in\Box$ means that the $(i,j)$-position in the matrix
    for $\widetilde{f}$ lies in the block $\Box$.

 \item [(2)]
  In summary$\,$:

  \fbox{\parbox{9em}{\flushleft\scriptsize
   \vspace{-3ex}
    Young diagram corresponding to a distinguished
    $S^1$-fixed-point component}}\
   $\Longrightarrow$
  \fbox{\parbox{9em}{\flushleft\scriptsize
   \vspace{-3ex}
    Collection of sparse lower-triangular block forms associated
    to the $S^1$-weight subspaces of a fiber of the normal bundle
    to the distinguished $S^1$-fixed-point component}}\
   $\Longrightarrow$
  \fbox{\parbox{9em}{\flushleft\scriptsize
   \vspace{-3ex}
   $P_0$-weight system $\Wt_{P_0}$ of each $S^1$-weight subspace}}
\end{itemize}

\bigskip

Let us now turn to the computation of
 the ${\Bbb C}^{\times}$-equivariant Euler class of the normal bundle
 to a distinguished component in $\Quot_{P(t)}({\cal E}^n)$.

\bigskip

\subsection{Structure of the induced bundle on
    $B{\largeBbb C}^{\times}
        \times_{{\scriptsizeBbb C}^{\times}}
         F_{\alpha_1,\,\ldots,\,\alpha_r\,;\,0,\,\ldots,\,0}$
    and the ${\largeBbb C}^{\times}$-euquvariant Euler class
     $e_{{\scriptsizeBbb C}^{\times}}\, 
       \nu_{F_{\alpha_1,\,\ldots,\,\alpha_r\,;\,0,\,\ldots,\,0}}
       \Quot_{P(t)}({\cal E}^n)$.                   }

We compute first the Chern polynomials of the normal bundle
 $\nu_{F_{\alpha_1,\,\ldots,\,\alpha_r\,;\,0,\,\ldots,\,0}}
                                    \Quot_{P(t)}({\cal E}^n)$
 of $F_{\alpha_1,\,\ldots,\,\alpha_r\,;\,0,\,\ldots,\,0}$
 in $\Quot_{P(t)}({\cal E}^n)$ and then
 use this to express the ${\Bbb C}^{\times}$-equivariant Euler class
 after working out the bundle structure of the induced bundle
 of $\nu_{F_{\alpha_1,\,\ldots,\,\alpha_r\,;\,0,\,\ldots,\,0}}$
 over
 $B{\Bbb C}^{\times}
   \times F_{\alpha_1,\,\ldots,\,\alpha_r\,;\,0,\,\ldots,\,0}$.

\bigskip

The following fact is in Borel and Hirzebruch [B-H], with slight
 modification to fit into our situation$\,$:

\bigskip

\noindent
{\bf Fact 3.3.1 [Chern class and representation].}
(Cf.\ [B-H], [B-T], [Fu1], [Hi], [M-S], and [Sp].) {\it 
\begin{itemize}
 \item [{\rm (1)}]
  Let
   $T$ be a maximal torus of $U(n)$,
   ${\frak h}$ be the corresponding Cartan subalgebra,
   and $\Fl(n):=U(n)/\mbox{\raisebox{-.4ex}{$T$}}$.
  Then, there are canonical homomorphisms
   $$
    \{\,\mbox{\rm integral linear functional on ${\frak h}$}\,\}\;
     \simeq\; H^1(T,{\Bbb Z})\; \rightarrow\; H^2(\Fl(n),{\Bbb Z})\,,
   $$
   where the second homomorphism is surjective and is given
   by the transgression homomorphism associated to the principal
   $T$-bundle $U(n)\rightarrow \Fl(n)$ from the quotient map.
  With respect to the defining representation of $U(n)$ on
   the Hermitian ${\Bbb C}^n$, $T$ corresponds to a unique orthonormal
   basis in ${\Bbb C}^n$ up to permutations. In terms of this basis,
   $T$ is realized as the group of unitary diagonal matrices. Thus,
   $T$ comes with a natural product decomposition $T=U(1)^{\times n}$
   that is invariant under the Weyl group action and each $U(1)$-factor
   of which is canonically oriented.
  This decomposition specifies then a distinguished basis
   $x_1,\,\ldots,\,x_n$ for $H^1(T,{\Bbb Z})$, unique up permutations.
  Regard $x_i$ also as elements in the other two groups
   via the above homomorphism and
   let $y_i=-x_i$ in $H^2(\Fl(n),{\Bbb Z})$.
  Up to permutations, $y_i$ in $H^2(\Fl(n),{\Bbb Z})$ are
   the first Chern class of the line bundles associated to
   the tautological flag bundle over $\Fl(n)$.
  These $y_i$ generate $H^2(\Fl(n),{\Bbb Z})$ and they satisfy
   $$
    \sigma_k(y_1,\,\ldots,\,y_n)\;=\; 0\,,\;
      \mbox{for $k\,=\,1,\,\ldots,\, n$}\,,
   $$
   where $\sigma_k$ is the elementary symmetric polynomial of degree
   $k$ for $n$ variables.

 \item [{\rm (2)}]
  {\rm [Chern root].}
  Let
   $P_0=U(m_1)\times\,\cdots\,\times U(m_k)\,\times U(m_{k+1})
                                                      \subset U(n)$,
    where \newline $m_1+\,\cdots\,+m_k+m_{k+1}=n$,
   $T$ be a maximal torus of $U(n)$ contained in $P_0$, and
   $\eta\,:\,U(n)\,
           \rightarrow\,B=U(n)/\mbox{\raisebox{-.4ex}{$P_0$}}$
   be the principal $P_0$-bundle over $B$ from the quotient map.
   Then $\Fl(n)$ is a split manifold for $\eta$.
  Let $\zeta: \Fl(n)\rightarrow B$ be the induced map from $\eta$, 
   then
   $\zeta^{\ast}:H^{\ast}(B,{\Bbb Z})
                   \rightarrow H^{\ast}(\Fl(n),{\Bbb Z})$
   is injective and $\zeta^{\ast}c(\eta)=\prod_{i=1}^n(1+y_i)$.

 \item [{\rm (3)}]
  {\rm [naturality of Chern class].}
  Let $\rho$ be an $m$-dimensional unitary representation of $P_0$
   with weights
   $w_j=a_{j1}x_1+\,\cdots\,+a_{jn}x_n$, $j=1,\,\ldots,\,m$, 
   and $V:= U(n)\times_{\rho}{\Bbb C}^m$ be the associated homogeneous
   vector bundle over $U(n)/\mbox{\raisebox{-.4ex}{$P_0$}}$.
  Then
   $$
    \zeta^{\ast}c(V)\;
     =\; \prod_{j=1}^m(1+w_j)\;
     =\; \prod_{j=1}^m(1+a_{j1}y_1+\,\cdots\,+a_{jn}y_n)\,.
   $$
 This expression is invariant under the Weyl group action $P_0$,
  $\Sym_{m_1}\times\,\cdots\,\times \Sym_{m_{k+1}}$ on the set
  $\{\, y_1,\,\ldots,\,y_{m_1};\:
        y_{m_1+1},\,\ldots,\, y_{m_1+m_2};\,\cdots;\:
        y_{m_1+\,\cdots\,+m_{k-1}+1},\,\ldots,\,y_{n-r}\,\}$
  by the permutations that
   $\Sym_{m_1}$ permutes the first $m_1$ letters, $\Sym_{m_2}$
   the next $m_2$ letters, and so on.
 The result is an integral polynomial function of symmetric functions
   in $y_1,\,\ldots,\,y_{m_1}$, in $y_{m_1+1},\,\ldots,\,y_{m_1+m_2}$,
   and so on respectively.
 Each of these partial symmetric products of Chern roots $y_i$ can be
  identified with the special Schubert cycles in the flag manifold
  $U(n)/\mbox{\raisebox{-.4ex}{$P_0$}}$.
\end{itemize}
} 

\bigskip

Recall the $P_0$-weight system associated to the Young diagram
 corresponding to $F_{\alpha_1,\,\ldots,\,\alpha_r\,;\,0,\,\ldots,\,0}$
 and Corollary 2.3.4, which says that all the $S^1$-weight of a fiber
 of the normal bundle are non-zero.
Let $\nu$ be the normal bundle in consideration.
Then the above fact implies that the Chern polynomial $c_{\nu}(t)$
 of $\nu$ is given by 
 $$
  c_{\nu}(t)\,
    =\;\prod_{
         \left(\,
           \begin{array}{c}
            \mbox{\footnotesize triangular block form}\:\Delta_w\\
            \mbox{\footnotesize for non-zero $S^1$-weight}
           \end{array}
         \right)}\;
       \prod_{(\,\mbox{\footnotesize block}\:\Box\in\Delta_w\,)}\;
       \prod_{(i,j)\in\Box}\,(t+y_{\,i+m_0}-y_j)\,,
 $$
 where the first product on the right hand side of the equality
  ranges over all possible non-zero $S^1$-weights $w$.
 The result is an integral polynomial function of
  the special Schubert cycles in
  $A_{\ast}(F_{\alpha_1,\,\ldots,\,\alpha_r\,;\,0,\,\ldots,\,0})$,
  cf.\ Fact/Definition 3.1.3.

The $S_1$-action on
 $\nu=\nu_{F_{\alpha_1,\,\ldots,\,\alpha_r\,;\,0,\,\ldots,\,0}}
   \Quot_{P(t)}({\cal E}^n)$
 induces a bundle
 $$
  {\cal T}\rightarrow
   \CP^{\infty}
    \times F_{\alpha_1,\,\ldots,\,\alpha_r\,;\,0,\,\ldots,\,0}\,.
 $$
Let $\nu=\oplus_w\,\nu_w$
 be the decomposition of the normal bundle as a direct sum of
 $S^1$-weight subspace and ${\cal T}=\oplus_w\,{\cal T}_w$
 be the induced decomposition of ${\cal T}$.

\bigskip

\noindent
{\bf Lemma 3.3.2 [induced bundle of $S^1$-weight summand].} {\it
  Let
   $$
    \CP^{\infty}\; \stackrel{{\rm pr}_1}{\longleftarrow}\;
    \CP^{\infty}
       \times F_{\alpha_1,\,\ldots,\,\alpha_r\,;\,0,\,\ldots,\,0}\;
      \stackrel{{\rm pr}_2}{\longrightarrow}\;
    F_{\alpha_1,\,\ldots,\,\alpha_r\,;\,0,\,\ldots,\,0}
   $$
   be the projection maps.
  Then
   $$
   {\cal T}_w\;
   =\; \pr_1^{\ast}\,
           {\cal O}_{{\scriptsizeBbb C}{\rm P}^{\infty}}(-w)
         \otimes\, \pr_2^{\ast}\,\nu_w\,.
   $$
} 

\bigskip

\noindent
{\it Proof.}
 Let
  $E=E{\Bbb C}^{\times}\rightarrow B{\Bbb C}^{\times}=\CP^{\infty}$
  be the universal principal ${\Bbb C}^{\times}$-bundle. 
 First notice that the associated line bundle of $E{\Bbb C}$
  to the representation of ${\Bbb C}^{\times}$ on ${\Bbb C}$
  by $v\mapsto tv$ for $t\in{\Bbb C}^{\times}$, $v\in{\Bbb C}$
  (i.e.\ the $w=1$ representation) is
  ${\cal O}_{\CPscriptsize^{\infty}}(-1)$.
 Since ${\Bbb C}^{\times}$ acts on $\nu_w$ by a single weight $w$,
  the induced action of ${\Bbb C}^{\times}$ on the projectivization
  ${\Bbb P}\nu_w$ of $\nu_w$ is trivial.
 Thus, as bundles over
  $\CP^{\infty}
    \times F_{\alpha_1,\,\ldots,\,\alpha_r\,;\,0,\,\ldots,\,0}$,
  $$
   {\Bbb P}{\cal T}_w\;
    =\; {\Bbb P}(E\times_{{\scriptsizeBbb C}^{\times}}\nu_w)\;
    =\; E\times_{{\scriptsizeBbb C}^{\times}}{\Bbb P}\nu_w\;
    =\; \CP^{\infty}\times {\Bbb P}\nu_w\;
    =\; {\Bbb P}\pr_2^{\ast}\nu_w\,.
  $$
 Since
  $1 \rightarrow {\Bbb C}^{\times} \rightarrow GL({\Bbb C})
     \rightarrow  \PGL({\Bbb C})   \rightarrow 1$
  is a central extension, the above isomorphism of projective
  bundles implies that
  ${\cal T}_w={\cal L}\otimes\pr_2^{\ast}\nu_w$
  for some line bundle ${\cal L}$ over \newline
  $\CP^{\infty}
   \times F_{\alpha_1,\,\ldots,\,\alpha_r\,;\,0,\,\ldots,\,0}$.
 By construction,
  $$
   {\cal T}_w|_{{\scriptsizeBbb C}{\rm P}^{\infty}\times\ast}\;
    \simeq\; \pr_1^{\ast}{\cal O}(-w)\otimes{\Bbb C}^R
    \hspace{1em}\mbox{and}\hspace{1em}
   {\cal T}_w
    |_{\ast\times
          F_{\alpha_1,\,\ldots,\,\alpha_r\,;\,0,\,\ldots,\,0}}\;
    \simeq\; \pr_2^{\ast}\nu_w\,,
  $$
  where $R$ is the rank of $\nu_w$.
 Since line bundles over flag manifolds are determined by their
   first Chern class and the second cohomology of flag manifolds are
   torsion-free,
  by the multiplicativity of the Chern character under tensor products
  and a comparison of first Chern classes, one concludes that 
  $$
   {\cal L}|_{{\scriptsizeBbb C}{\rm P}^{\infty}\times\ast}\;
    \simeq\;{\cal O}_{{\scriptsizeBbb C}{\rm P}^{\infty}}(-w)
   \hspace{1em}\mbox{and}\hspace{1em}
   {\cal L}
    |_{\ast\times F_{\alpha_1,\,\ldots,\,\alpha_r\,;\,0,\,\ldots,\,0}}\;
     \simeq\;
    {\cal O}_{F_{\alpha_1,\,\ldots,\,\alpha_r\,;\,0,\,\ldots,\,0}}\,.
  $$

 Consider now a finite model $\CP^N$ for $\CP^{\infty}$ with $N$
  very large.
 Then, since
  $$
   X_N \;
    :=\; \CP^N
         \times F_{\alpha_1,\,\ldots,\,\alpha_r\,;\,0,\,\ldots,\,0}
  $$
  is simply-connected and K\"{a}hler, from the long exact sequence
  {\scriptsize
  $$
   \hspace{-6em}
   \begin{array}{ccccccccccc}
    \cdots & \longrightarrow  & H^1(X_N,{\scriptsizeBbb Z})
           & \longrightarrow  & H^1(X_N, {\cal O}_{X_N})
           & \longrightarrow  & H^1(X_N, {\cal O}^{\ast}_{X_N})
           & \stackrel{c_1}{\longrightarrow}
           & H^2(X_N,{\scriptsizeBbb Z})
           & \longrightarrow  & \cdots \\[.6ex]
    &      & \|   & & \|    & &   \|    & &   \|  & &   \\[.6ex]
    &      & 0    & &  H^{0,1}_{\overline{\partial}}(X_N)=0
           &      & \Pic(X_N) &
           &  H^2(\CPscriptsize^N,{\scriptsizeBbb Z})
              \oplus
              H^2(
               F_{\alpha_1,\,\ldots,\,\alpha_r\,;\,0,\,\ldots,\,0},
               {\scriptsizeBbb Z})  &   &
   \end{array}
  $$
  {\normalsize associated}}  
  to the exponential sequence
  $0 \rightarrow {\Bbb Z} \rightarrow {\cal O}_{X_N}
     \rightarrow {\cal O}^{\ast}_{X_N} \rightarrow 0$,
  one concludes that the Picard variety $\Pic(X_N)$ is contained in
  $\Pic(\CP^N)\times
         \Pic(F_{\alpha_1,\,\ldots,\,\alpha_r\,;\,0,\,\ldots,\,0})$
  and hence that every line bundle on $X_N$ is of the form
  $\pr_1^{\ast}{\cal L}_1\otimes\pr_2^{\ast}{\cal L}_2$.
 Together with the earlier discussion in the proof, one has
  in particular that
  ${\cal L}
   =\pr_1^{\ast}\,{\cal O}_{{\scriptsizeBbb C}{\rm P}^N}(-w)$
  over $X_N$ for all large $N$.
 Let $N\rightarrow\infty$, one then concludes the lemma.

\noindent\hspace{14cm} $\Box$

\bigskip

Let $R(w)$ be the rank of $\nu_w$. By the multiplicativity of
 Euler class and the rule under the tensor with a line bundle
 (cf.\ [Fu1]), we conclude that
\bigskip

\noindent 
{\bf Theorem 3.3.3 [Euler class].} {\it
The $S^1$ equivariant Euler class of the normal bundle
 $\nu=\nu_{F_{\alpha_1,\,\ldots,\,\alpha_r\,;\,0,\,\ldots,\,0}}
   \Quot_{P(t)}({\cal E}^n)$
is given by
 \begin{eqnarray*}
  \lefteqn{ 
    e_{{\tinyBbb C}^{\times}}\,\nu\;
     =\; \prod_w\, e_{{\tinyBbb C}^{\times}}\nu_w\;
     =\; \prod_w\,
           c_{\nu_w}\,(-w\alpha)                } \\[.6ex]
    & & =\;\prod_{
             \left(\,
              \begin{array}{c}
               \mbox{\footnotesize triangular block form}\:\Delta_w\\
               \mbox{\footnotesize associated to nonzero $w$}
              \end{array}
             \right)}\;
          \prod_{(\,\mbox{\footnotesize block}\:\Box\in\Delta_w\,)}\;
          \prod_{(i,j)\in\Box}\,(\,-w\alpha\,+\,y_{i+m_0}\,-\,y_j\,)\,,
  \end{eqnarray*}
 where
  $\alpha=c_1({\cal O}_{\CPscriptsize^{\infty}}(1))$,
  $c_{\nu_w}(t)=t^{R(w)}+c_1(\nu_w)t^{R(w)-1}+\cdots\;$
  is the Chern polynomial of $\nu_w$. 
 The triangular block forms $\Delta_w$,
 $-\alpha_r\leq w\leq max\{\alpha_r-\alpha_1-1,0\}$,
 associated with the Young diagram corresponding to
 the partition $d=\alpha_1+\cdots+\alpha_r$
 are defined by Item {\rm (4)} in the subheading 
 ``$S^1$-weight-subspace decomposition" of the heading
 ``The $P_0$-module decomposition of $S^1$-weight spaces and
   the $P_0$-weight system" in Sec.\ 3.2. 
} 

\bigskip
Here  $\pr_1^{\ast}$, $\pr_2^{\ast}$ in the formula are omitted for
   simplicity of notations.

\bigskip

\section{Illustrations by two examples.}

In this section, we present two simple examples of
 the Mirror Principle computation that are computable by hand to
 illustrate the discussions in this article.
In these examples, the distinguished $S^1$-fixed-point components
 in the related components of Quot-schemes are either Grassmannian
 manifolds or complete flag manifolds.
The Schubert calculus of these follow from Fulton in [Fu1] and
 Monk in [Mo].
In particular, for the complete flag manfold
 $\Fl(3):=\Fl_{1,2}({\Bbb C}^3)$, the cohomology ring
 $H^{\ast}(\Fl(3),{\Bbb Z})$ is generated by $y_1$, $y_2$,
 ($y_3=-(y_1+y_2)$), where $y_i$ are the first Chern class the graded
 line bundles on $\Fl(3)$ associated to the flag of universal rank-$1$
 and rank-$2$ subbundles over $\Fl(3)$, cf.\ Fact 3.3.1.
The intergral of the top classes are given by
 $$
  \int_{Fl(3)}\, y_1^3\;=\; \int_{Fl(3)}\,y_2^3\;=\;0
   \hspace{1em}\mbox{and}\hspace{1em}
  \int_{Fl(3)}\, y_1^2y_2\;=\; -\,\int_{Fl(3)}\, y_1y_2^2\;=\; -1\,,
 $$
following [Mo].

\bigskip

\noindent
{\bf Example 4.1 [$\Gr_2({\Bbb C}^3)$, degree 3].}
In this case, $n=3$, $r=2$, $n-r=1$, and $d=3$.
There are two distinguished $S^1$-fixed-point components
 in the related component $\Quot_{P(t)=t+4}({\cal E}^3)$
 of Quot-scheme:

\begin{itemize}
 \item [$\bullet$]
  $F_{\,0,\,3\;;\;0,0}\, \simeq\,\Fl_{1,\,2}({\smallBbb C}^3)$,
  $\dimm\,=\,3$.
\end{itemize}

\centerline{\psfig{figure=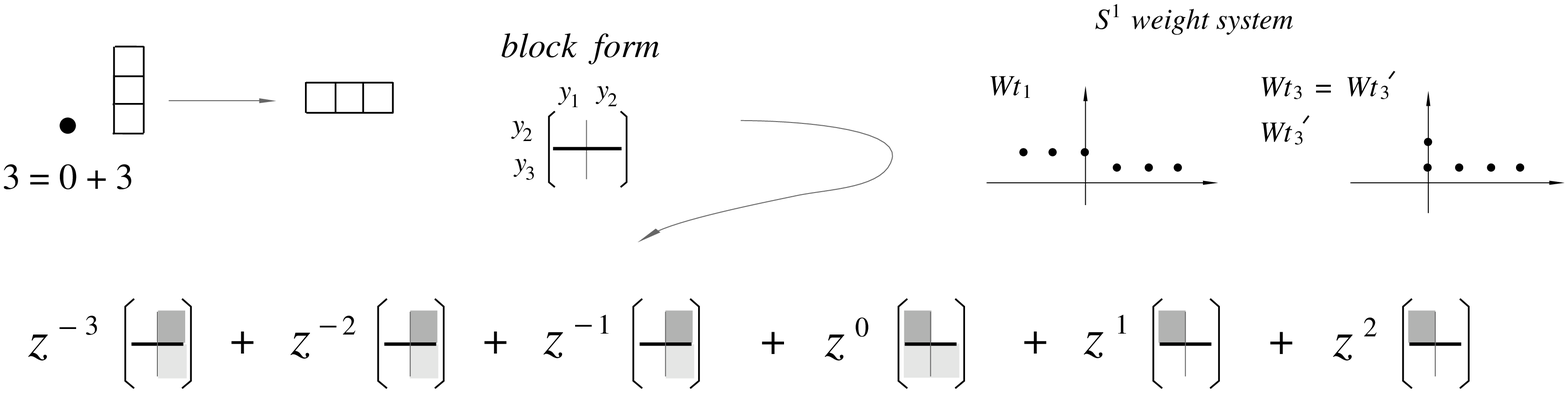,width=15cm,caption=}}

\begin{itemize}
 \item [${\boldmath -}$]
  Grouping of Chern roots$\,$: $\{\,y_1\;;\;y_2\;;\; y_3\,\}$.

 \item [${\boldmath -}$]
  ${\Bbb C}^{\times}$-equivariant Euler class of normal bundle$\,$:
  \begin{eqnarray*}
   \lefteqn{e_{{\tinyBbb C}^{\times}}(\nu)\;
    =\; (-3\alpha)\,(-3\alpha+y_3-y_2)\,
        (-2\alpha)\,(-2\alpha+y_3-y_2)} \\
    & &  \hspace{6em} \cdot\,
          (-\alpha)\,(-\alpha+y_3-y_2)\,
          (\alpha+y_2-y_1)\,(2\alpha+y_2-y_1)\,.
  \end{eqnarray*}

 \item [${\boldmath -}$]
  Pulled-back hyperplane class$\,$:
  $$
   k^{\ast}\psi^{\ast}\kappa\;=\; g^{\ast}j^{\ast}\kappa\;
   =\;-(y_1+y_2)\,.
  $$

 \item [${\boldmath -}$]
  The integral over the component$\,$:
  $$
   \int_{E}\,\frac{k^{\ast}\psi^{\ast}e^{\kappa\cdot\zeta}}{
               e_{{\tinyBbb C}^{\times}}(E/{\cal Q}_d)}\;
   =\; -\, \frac{103}{1296}\,\frac{1}{\alpha^{11}}\,
       -\,\frac{23}{108}\,\frac{\zeta}{\alpha^{10}}\,
       -\frac{29}{864}\,\frac{\zeta^2}{\alpha^9}\,.
  $$
\end{itemize}

\begin{itemize}
 \item [$\bullet$]
  $F_{\,1,\,2\;;\;0,0}\, \simeq\,\Fl_{1,\,2}({\smallBbb C}^3)$,
  $\dimm\,=\,3$.
\end{itemize}

\centerline{\psfig{figure=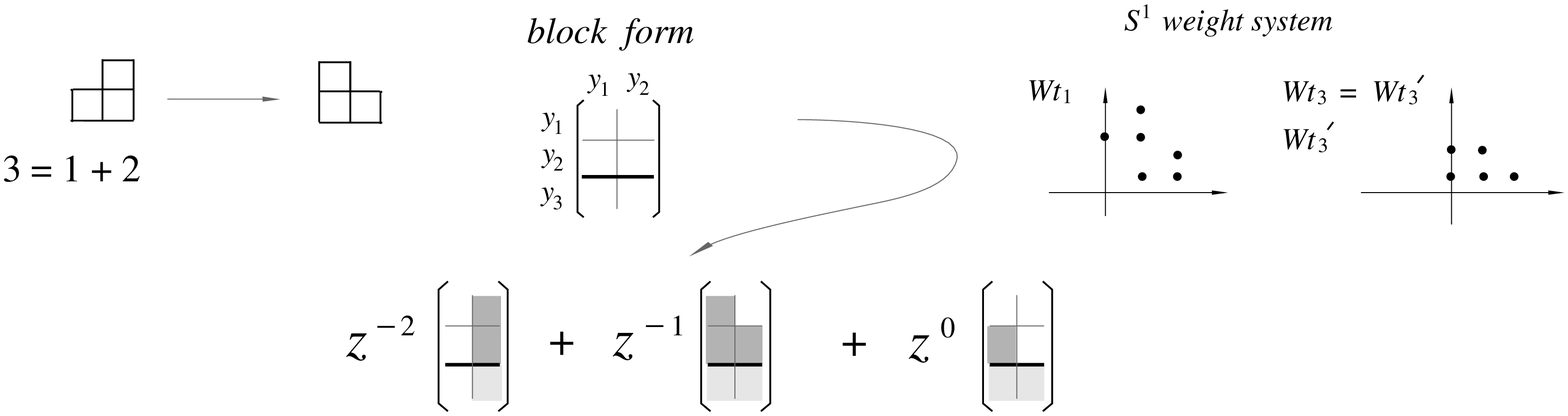,width=15cm,caption=}}

\begin{itemize}
 \item [${\boldmath -}$]
  Grouping of Chern roots$\,$: $\{\,y_1\;;\;y_2\;;\; y_3\,\}$.

 \item [${\boldmath -}$]
  ${\Bbb C}^{\times}$-equivariant Euler class of normal bundle$\,$:
  \begin{eqnarray*}
   \lefteqn{e_{{\tinyBbb C}^{\times}}(\nu)\;
      =\; (-2\alpha+y_1-y_2)\,(-2\alpha)\,(-2\alpha+y_3-y_2)} \\
     & & \hspace{6em}\cdot\,
       (-\alpha)\,(-\alpha+y_2-y_1)\,(-\alpha+y_3-y_1)\,(-\alpha)\,
        (-\alpha+y_3-y_2)\,.
  \end{eqnarray*}

 \item [${\boldmath -}$]
  Pulled-back hyperplane class$\,$:
  $$
   k^{\ast}\psi^{\ast}\kappa\;=\; g^{\ast}j^{\ast}\kappa\;
   =\;-(y_1+y_2)\,.
  $$

 \item [${\boldmath -}$]
  The integral over the component$\,$:
  $$
   \int_{E}\,\frac{k^{\ast}\psi^{\ast}e^{\kappa\cdot\zeta}}{
               e_{{\tinyBbb C}^{\times}}(E/{\cal Q}_d)}\;
   =\;\frac{3}{16}\,\frac{\zeta}{\alpha^{10}}\,
      +\,\frac{1}{32}\,\frac{\zeta^2}{\alpha^9}\,.
  $$
\end{itemize}

$$
 \mbox{Total integral}\;
   =\; -\,\frac{103}{1296}\,\frac{1}{\alpha^{11}}\,
       -\,\frac{11}{432}\,\frac{\zeta}{\alpha^{10}}\,
       -\frac{1}{432}\,\frac{\zeta^2}{\alpha^9}\,.
$$

\noindent\hspace{14cm}$\Box$

\bigskip

\noindent
{\bf Example 4.2 [$\Gr_1({\Bbb C}^3)$, degree 3].}
In this case, $n=3$, $r=1$, $n-r=2$, and $d=3$.
There is one distinguished $S^1$-fixed-point component in
 the related component $\Quot_{P(t)=2t+5}({\cal E}^3)$ of Quot-scheme:

\begin{itemize}
 \item [${\boldmath -}$]
  $F_{\,3\;;\;0}\, \simeq\,\Gr_1({\smallBbb C}^3)$, $\dimm\,=\,2$.
\end{itemize}

\centerline{\psfig{figure=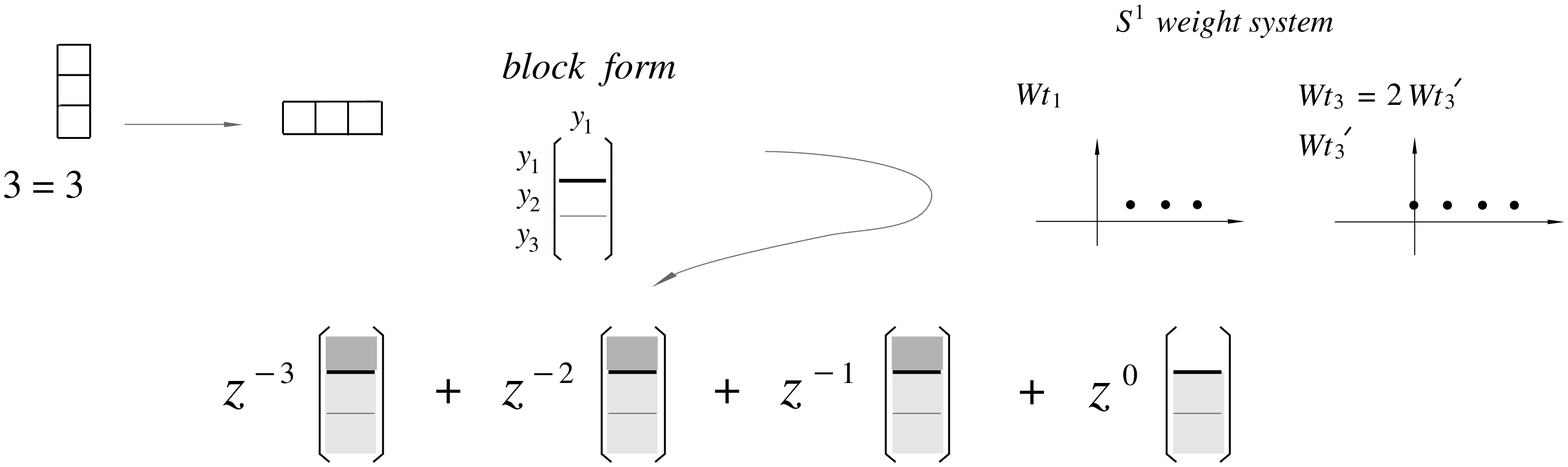,width=15cm,caption=}}

\begin{itemize}
 \item [${\boldmath -}$]
  Grouping of Chern roots$\,$: $\{\,y_1\;;\;y_2, y_3\,\}$.

 \item [${\boldmath -}$]
  ${\Bbb C}^{\times}$-equivariant Euler class of normal bundle$\,$:
  \begin{eqnarray*}
   \lefteqn{e_{{\tinyBbb C}^{\times}}(\nu)\;
     =\; (-3\alpha)\,(-3\alpha+y_2-y_1)\,(-3\alpha+y_3-y_1) 
         (-2\alpha)\,(-2\alpha+y_2-y_1)  } \\
        & & \hspace{8em}\cdot\,
            (-2\alpha+y_3-y_1)\,
            (-\alpha)\,(-\alpha+y_2-y_1)\,(-\alpha+y_3-y_1)\,.
  \end{eqnarray*}

 \item [${\boldmath -}$]
  Pulled-back hyperplane class$\,$:
  $$
   k^{\ast}\psi^{\ast}\kappa\;=\; g^{\ast}j^{\ast}\kappa\;=\;-y_1\,.
  $$

 \item [${\boldmath -}$]
  The integral over the component$\,$:
  $$
   \int_{E}\,\frac{k^{\ast}\psi^{\ast}e^{\kappa\cdot\zeta}}{
               e_{{\tinyBbb C}^{\times}}(E/{\cal Q}_d)}\;
   =\; -\,\frac{103}{1296}\,\frac{1}{\alpha^{11}}\,
       -\,\frac{11}{432}\,\frac{\zeta}{\alpha^{10}}\,
       -\frac{1}{432}\,\frac{\zeta^2}{\alpha^9}\,,
  $$
  which is the same as the total integral in Example 4.1,
  as it should be since \newline
  $\Gr_1({\Bbb C}^3)=\Gr_2({\Bbb C}^3)$.
\end{itemize}

\noindent\hspace{14cm}$\Box$

\bigskip

\noindent
{\it Remark 4.3.}
One can check that the integral values are correct,
 using the result in [L-L-Y1, I] for the computation for $\CP^2$.
Simple examples as they are, one observes that the intermediate
 details in the computation do depend on the presentation of
 a Grassmannian manifold and these details are in general
 very different.
The fact that either presentation gives an identical answer
 provides thus a computational check of the theory developed.

\bigskip

\noindent
{\it Remark 4.4.}
Now that we can compute the integral that is related to the intersection
 numbers on the moduli space of rational stable maps into Grassmannian
 manifolds, the A-model for Calabi-Yau complete intersections 
 in a Grassmannian manifold can also be computed explicitly.

\bigskip

\noindent
{\it Remark 4.5.}
Two immediate questions follow from the current work$\,$:
\begin{itemize}
 \item[(1)]
   the automatization of the calculations via a computer
   code, following the diagrammatic rules discussed,
   and the computation for more examples and

 \item[(2)]
  generalization of the discussion to flag manifolds,
  which involves hyper-Quot schemes.
\end{itemize}
The study of them will be reported in another work.

\bigskip

With these remarks, we conclude this paper.

\newpage

\end{document}